\documentclass[11pt,leqno]{amsart}
\usepackage[T1]{fontenc}
\usepackage{lmodern}
\usepackage{amssymb,mathtools}
\usepackage{microtype}
\usepackage{xcolor}
\usepackage[hidelinks,bookmarksnumbered]{hyperref}
\hypersetup{pdftitle={The invariant subspace problem for trace-class perturbations of
subnormal operators}}
\hypersetup{pdfauthor={Kunyu Guo, Chunlan Jiang, Penghui Wang, Zeyou Zhu}}
\newtheorem{theorem}{Theorem}[section]
\newtheorem{proposition}[theorem]{Proposition}
\newtheorem{lemma}[theorem]{Lemma}

\newtheorem{Problem}[theorem]{Problem}

\numberwithin{equation}{section}
\title[Invariant subspace problem for trace-class perturbations of subnormal operators]
{Invariant subspace problem for trace-class\protect\\
 perturbations of subnormal operators}
\author{Kunyu Guo}
\address{Kunyu Guo: School of Mathematical Sciences, Fudan University, Shanghai, 200433, China}
\email{kyguo@fudan.edu.cn}

\author{Chunlan Jiang}
\address{Chunlan Jiang: School of Mathematical Sciences, Hebei Normal University, Shijiazhuang, 050024, China}
\email{cljiang@hebtu.edu.cn}

\author{Penghui Wang}
\address{Penghui Wang: School of Mathematics, Shandong University, Jinan, 250100, China}
\email{phwang@sdu.edu.cn}

\author{Zeyou Zhu}
\address{Zeyou Zhu: School of Mathematical Sciences, Fudan University, Shanghai, 200433, China}
\email{zhuzeyou@fudan.edu.cn}

\date{}

\begin{document}
\begin{abstract}
Let $S$ be a subnormal operator on a separable infinite-dimensional complex
Hilbert space, and let $K$ be a trace class operator. In the present paper, we prove that  $S+K$ has a nontrivial invariant subspace. Furthermore, for a normal operator $S$, we prove that if the Hausdorff dimension of $\sigma(S)$ is greater than $1$, then $S+K$ has a hyperinvariant half-space.
\end{abstract}
\maketitle


\section{Introduction}\label{sec:intro}

Let $H$ be a separable infinite-dimensional complex Hilbert space.
As a special case of the invariant subspace problem,
Halmos~\cite[p.~899]{Halmos1970} asked the following question.
\begin{Problem}\label{intro:normal-rank-one}
If $N$ is normal and $K$ has rank one, is $N+K$ intransitive?
\end{Problem}
An operator $T\in\mathcal B(H)$ is \emph{intransitive} if there is a
closed subspace $\{0\}\subsetneq\mathcal M\subsetneq H$ such that
$T\mathcal M\subseteq\mathcal M$.
Pearcy~\cite[Chapter~8, Problem~K]{Pearcy} asked whether $D+F$ is
intransitive for every rank-one operator $F$ when $D$ is diagonal and
its diagonal entries form a dense subset of the unit disk.
For \textcolor{black}{work on} this problem \textcolor{black}{and related questions}, see \cite{AlbrechtChevreau2021,Berg1971,Chalendar1997,Donoghue,FX,FJKP,FoiasJungKoPearcy2008II,FoiasJungKoPearcy2011,
GallardoGonzalez2022,GallardoGonzalez2023idempotent,GGD24,GGD25,Ionascu2001,JungLeeSeo2018,Kitano1968,Klaja2015rankone,
Klaja2015hyperinvmult,PY,Radjabalipour1975,RadjabalipourRadjavi1975,Schwartz1962,Stampfli1984} and the references therein.
Recently, the second author and his collaborators have developed several
approaches to the invariant subspace problem~\cite{FangHouJiang2024,FHJZ1,FHJZ,JiangFangJi2025,ChengFangJiang2025,JiJiangWang2026,JiangCao2023}.
In Pearcy's diagonal case of Problem~\ref{intro:normal-rank-one},
let \(T=D+f\otimes g\), where \(D=\operatorname{diag}(\lambda_0,\lambda_1,\ldots)\),
, as a consequence,
Fang, Hou, Jiang, and Zhang~\cite{FHJZ} proved that \(T\) is intransitive if either
\(f\) or \(g\) is noncyclic for \(D\).

The spectral theorem implies that every normal operator on $H$ is
intransitive. Brown~\cite{Brown1978} proved that every subnormal
operator on $H$ is intransitive. Here we consider trace-class
perturbations of subnormal operators. Write $S_p(H)$ for the
Schatten-$p$ class on $H$. In particular, $S_1(H)$ is the trace class
and $S_2(H)$ is the Hilbert--Schmidt class. Our first main result is
the following.

\begin{theorem}\label{main:invariant}
Let $S$ be a subnormal operator on $H$ and let $K\in S_1(H)$.
Then $S+K$ is intransitive.
\end{theorem}

Since every normal operator is subnormal and every rank-one operator is
trace class, Theorem~\ref{main:invariant} answers Halmos's rank-one
question affirmatively and, in particular, resolves Pearcy's problem.

Our second result concerns hyperinvariant subspaces of $N+K$ when
$N$ is normal. A closed subspace $\mathcal L$ is \emph{hyperinvariant}
for $T$ if $S\mathcal L\subseteq\mathcal L$ for every $S\in\{T\}'$.
It is a \emph{half-space} if
$\dim\mathcal L=\dim\mathcal L^\perp=\infty$.
% To state our second main result, let $\mu$ be a probability measure on $\mathbb C$ with compact supports and write the Riesz $1$-capacity of $\mu$ as
%\[
% I_1(\mu)=\iint\frac{d\mu(z)\,d\mu(w)}{|z-w|}.
%\]
%A compact set has positive Riesz $1$-capacity
%if it supports a probability measure of finite $1$-energy.

\begin{theorem}\label{main:half}
Let $N$ be normal and let $K\in S_1(H)$. If
$\dim_{\mathrm H}(\sigma(N))>1$, then $N+K$ has a hyperinvariant
half-space, where $\dim_{\mathrm H}$ denotes Hausdorff dimension.
\end{theorem}
In particular, if the diagonal entries of $D$ form a dense subset of
the unit disk, then $\sigma(D)=\overline{\mathbb D}$ and
$\dim_{\mathrm H}(\sigma(D))=2$. Thus, in Pearcy's setting, $D+K$
has a hyperinvariant half-space for every $K\in S_1(H)$.

%On the other hand, a countably infinite spectrum
%can occur for a normal $N$ and a trace-class $K$ such that every
%nonzero closed invariant subspace of $N+K$ has finite codimension
%(Section~\ref{sec:examples}).

%To relate the two theorems, put $T=N+K$, $\Lambda=\sigma(N)$ and
%$\|p\|_L=\sup_{z\in L}|p(z)|$. Suppose that a bounded disk $U$ contains
%a compact subset of $\Lambda$ of positive analytic capacity.
%Proposition~\ref{gen:local} produces $x\ne0$ for which
%\begin{equation}\label{intro:local}
% \|p(T)x\|\le C_x\|p\|_U\qquad (p\in\mathbb C[z]).
%\end{equation}
%The construction uses the Cauchy transform of $\mu=BE(\,\cdot\,)A$,
%where $K=AB$ with $A,B$ Hilbert--Schmidt and $E$ is the
%spectral measure of $N$. Relative to a positive measure $\nu$
%supported in $U\cap\Lambda$, write $\mu=\mathcal M\,\nu+\mu_s$ and
%$\Delta=I-\mathcal C\mu$. Where $\Delta$ is invertible and $\mathcal M\ne0$, the choice
%$f=\Delta^{-1}\mathcal M e$ for a suitable basis vector $e$ supplies the correction term in
%\eqref{gen:correction}. Disjoint disks for $T$ and $T^*$ then give
%orthogonal families of vectors. This proves the half-space conclusion
%under the more general assumption $\gamma(\Lambda)>0$
%(Proposition~\ref{eng:positive}). Under the finite-energy hypothesis
%the local vectors also admit a construction using absolutely convergent
%integrals.

This paper is organized as follows. Section~\ref{sec:cauchy} gives
estimations for Cauchy transforms in terms of analytic capacity. In
Section~\ref{sec:model}, we use spectral integrals to obtain local
polynomial bounds for $T=N+K$. In Section~\ref{sec:localization},
we obtain these bounds for a nonzero vector on every disk that meets
$\sigma(N)$ in a set of positive analytic capacity.
Section~\ref{sec:half} uses these results to prove the main theorems on
intransitivity and hyperinvariant half-spaces.

\section{Cauchy transforms and analytic capacity}\label{sec:cauchy}

In this section, we introduce the Cauchy transforms for measure and give some estimations for analytic capacity. Let \(m_2\) be the Lebesgue measure on the complex plane \(\mathbb C\). For a separable complex Hilbert space $V$, let $\mu$ be a $V$-valued  Borel measure on $\mathbb C$ of
finite variation, where the total variation measure for $\mu$ is defined by
$$
|\mu|(A)=\sup\left\{\sum\limits_{j=1}^n\|\mu(A_j)\|_V;\, A=\sqcup_{j=1}^n A_j, \, A_j \text{ are Borel set}, \, n\geq1 \right\}; 
$$
moreover, we assume  $\mathrm{supp}(\mu)$ is a compact subset of $\mathbb C$.

For $\varepsilon>0$, let
\[
  \mathcal C_\varepsilon\mu(z)
  =\int_{|z-w|>\varepsilon}\frac{d\mu(w)}{z-w}.
\]
Whenever the norm limit exists, write
\[
  \mathcal C\mu(z)=\lim\limits_{\varepsilon\to0^+}\mathcal C_\varepsilon\mu(z).
\]

For any $z\in\mathbb C$, define the map $F^z: (0,\infty)\to V$  by $F^z(\varepsilon)=\mathcal C_\varepsilon\mu(z)$. 
The following lemma deals with the right continuity of $F^z$.
\begin{lemma}\label{cau:rational-truncations}
For each $z\in\mathbb C$, the map $F^z: (0,\infty)\to V$ is right-continuous, that is, for any $\varepsilon_0>0$, one has
$$
\lim\limits_{h\to 0^+} \|\mathcal C_{\varepsilon_0+h}\mu(z)-\mathcal C_{\varepsilon_0}\mu(z)\|_V=0.
$$
Consequently, for the maximal Cauchy transformation $\mathcal C_*\mu(z)
  :=\sup_{\varepsilon>0}\|\mathcal C_\varepsilon\mu(z)\|_V$, one has
\[
  \mathcal C_*\mu(z)
   =\sup_{q\in\mathbb Q\cap(0,\infty)}\|\mathcal C_q\mu(z)\|_V,
\]
which is
\([0,\infty]\)-valued Borel function.
\end{lemma}

\begin{proof}
For any fixed $z\in\mathbb C$ and $\varepsilon_0>0$. For $q>\varepsilon_0$,
\[
\begin{aligned}
  \limsup\limits_{q\to\varepsilon_0^+}\|F^z(q)-F^z(\varepsilon_0)\|_v
  &=\limsup\limits_{q\to\varepsilon_0^+}\|\mathcal C_{\varepsilon_0}\mu(z)-\mathcal C_q\mu(z)\|_V\\
  &\le  \limsup\limits_{q\to\varepsilon_0^+} \int_{\varepsilon_0<|z-w|\le q}
          \frac{d|\mu|(w)}{|z-w|} \\
  &\le \limsup\limits_{q\to\varepsilon_0^+} \varepsilon_0^{-1}
       |\mu|\{w:\varepsilon_0<|z-w|\le q\}=0.
\end{aligned}
\]
 For convenience, set \(\mathbb Q_+=\mathbb Q\cap(0,\infty)\). By right continuity of $F^z$,
\[
  \|\mathcal C_\varepsilon\mu(z)\|_V
  =\lim\limits_{\substack{q\to\varepsilon^+\\q\in\mathbb Q_+}}
    \|\mathcal C_q\mu(z)\|_V
  \le\sup_{q\in\mathbb Q_+}\|\mathcal C_q\mu(z)\|_V.
\]
Taking the supremum over \(\varepsilon>0\) and using
\(\mathbb Q_+\subset(0,\infty)\), we obtain
\[
  \sup_{\varepsilon>0}\|F^z(\varepsilon)\|_V=\sup_{\varepsilon>0}\|\mathcal C_\varepsilon\mu(z)\|_V
  =\sup_{q\in\mathbb Q_+}\|\mathcal C_q\mu(z)\|_V.
\]
For any fixed \(q>0\), set
\[
  \kappa_q(z,w)=
  \begin{cases}
    (z-w)^{-1},&|z-w|>q,\\
    0,&|z-w|\le q.
  \end{cases}
\]
For convenience, for any fixed $q>0$ write $F_q(z)=F^z(q)$. Then  
\[
  F_q(z)
  =\int_{\mathbb C}\kappa_q(z,w)\,d\mu(w)
\]
Since \(\kappa_q\) is a  Borel function and \(\|\kappa_q\|_\infty\le q^{-1}\), the map $F_q(\cdot)$
is a $V$-valued Borel function on $\mathbb C$. Thus, for every $t>0$,
\[
  \{\mathcal C_*\mu>t\}
  =\bigcup_{q\in\mathbb Q_+}\{z:\|\mathcal C_q\mu(z)\|_V>t\},
\]
and hence $\mathcal C_*\mu$ is a Borel function.
\end{proof}

Write \(\widehat{\mathbb C}=\mathbb C\cup\{\infty\}\) for the
Riemann sphere. For a compact set \(L\subset\mathbb C\), let
\[
  \gamma(L)=\sup\bigl\{
    |f'(\infty)| : f\in H^\infty(\widehat{\mathbb C}\setminus L),
    f(\infty)=0,\ \|f\|_\infty\le1
  \bigr\},
\]
where  \(H^\infty\) denotes bounded holomorphic functions and $f'(\infty)$ can be understood as follows.  Using the local coordinate $w=1/z$ near $\infty$, and consider $g(w)=f(1/w)$, and $g(0)=f(\infty)$. Since $f(\infty)=0$, we have 
$$
f'(\infty)=g'(0)=\lim_{w\to0} \frac{g(w)}{w}=\lim\limits_{z\to\infty}z f(z).
$$  For arbitrary
\(E\subset\mathbb C\), set
\[
  \gamma(E)=\sup\{\gamma(L):L\subset E,\ L\text{ compact}\}.
\]

For every disk $U=B(a,r)$,
\begin{equation}\label{cau:area}
\begin{aligned}
  c_U
  &:=\sup_{w\in\mathbb C}\int_U\frac{dm_2(z)}{|z-w|}\\
  &\le\sup_{w\in\mathbb C}
      \left(\int_{B(w,r)}\frac{dm_2(z)}{|z-w|}
            +\frac{m_2(U)}{r}\right)
   =2\pi r+\pi r=3\pi r.
\end{aligned}
\end{equation}
\begin{lemma}\label{cau:area-lemma}
Let \(\mu\) be a \(V\)-valued Borel measure of finite variation, and let
\[
  I_\mu(z)=\int_{\mathbb C}\frac{d|\mu|(w)}{|z-w|},
  \, \text{ where } \frac10=+\infty.
\]
For every open disk \(U\),
\[
  \int_U I_\mu\,dm_2\le c_U|\mu|(\mathbb C).
\]
Moreover,
\[
  \mathcal C\mu(z)
  =\lim\limits_{\varepsilon\to0^+}\mathcal C_\varepsilon\mu(z)
  =\int_{\mathbb C}\frac{d\mu(w)}{z-w}
  \qquad\text{for \(m_2\)-almost every }z\in\mathbb C,
\]
where the last integral converges absolutely.
\end{lemma}

\begin{proof}
By \eqref{cau:area} and Tonelli's theorem,
\[
\begin{aligned}
  \int_U I_\mu(z)\,dm_2(z)
  =\int_{\mathbb C}
    \left(\int_U\frac{dm_2(z)}{|z-w|}\right)d|\mu|(w)
  \le c_U|\mu|(\mathbb C)<\infty.
\end{aligned}
\]
Hence \(I_\mu\in L^1_{\mathrm{loc}}(\mathbb C,m_2)\), so
\(I_\mu(z)<\infty\) for \(m_2\)-almost every \(z\).
At each such point, \(|\mu|(\{z\})=0\) and
\[
  \int_{\mathbb C}\frac{d|\mu|(w)}{|z-w|}<\infty.
\]
Thus the Cauchy integral converges absolutely, and dominated
convergence gives
\[
 \limsup\limits_{\varepsilon\to0^+} \left\|\mathcal C_\varepsilon\mu(z)
       -\int_{\mathbb C}\frac{d\mu(w)}{z-w}\right\|_V
  \le\limsup\limits_{\varepsilon\to0^+}\int_{0<|z-w|\le\varepsilon}
       \frac{d|\mu|(w)}{|z-w|}=0.
\]
\end{proof}

For a finite complex Borel measure \(\xi\) on \(\mathbb C\), define its
truncated and maximal Cauchy transforms by
$$
\mathcal C_\varepsilon \xi (z)=\int_{|z-w|>\varepsilon} \frac{d\xi(w)}{z-w},
\text{ and }\mathcal C_*\xi(z)
  :=\sup_{\varepsilon>0}|\mathcal C_\varepsilon\xi(z)|.$$ 
  %The following lemma was obtained by Mateu, Tolsa and Verdera \cite[Section 2]{MTV} and \cite[Section 1.2]{CMPT}. 
Part~\textup{(1)} of the following lemma is recorded in
\cite[Section~2, following (2.7)]{MTV}, and part~\textup{(2)} follows
from the characterization in \cite[Section~1.2, (8)--(10)]{CMPT}.
\begin{lemma}\label{cau:capacity}
There exist constants \(c,C>0\), independent of the measure \(\xi\),
the sets \(L,F\), and the parameters \(t,z,r\), with the following
properties.
\begin{enumerate}
\item
For every finite complex Borel measure \(\xi\) on \(\mathbb C\),
\begin{equation}\label{cau:scalar}
  \gamma\{z:\mathcal C_*\xi(z)>t\}
  \le \frac{C|\xi|(\mathbb C)}{t}, 
  \text{ for } t>0.
\end{equation}

\item
For every compact \(L\subset\mathbb C\) with \(\gamma(L)>0\),
there is a finite positive measure \(\nu\) supported on \(L\)
such that
\begin{equation}\label{cau:measure}
  \nu(L)\ge c\gamma(L),\,
  \nu(B(z,r))\le r,\text{ and }
  \nu(F)\le C\gamma(F)
\end{equation}
for all \(z\in\mathbb C\), \(r>0\), and Borel sets
\(F\subset\mathbb C\).
\end{enumerate}
\end{lemma}
We explain only the third inequality in \eqref{cau:measure}.
For a compact set \(K\subset\mathbb C\), a finite positive Borel
measure \(\rho\) supported on \(K\) is called admissible if
\(\rho(B(z,r))\le r\) for all \(z\in\mathbb C\) and \(r>0\), and
\[
  \sup_{\varepsilon>0}
  \|\mathcal C_\varepsilon(f\rho)\|_{L^2(\rho)}
  \le \|f\|_{L^2(\rho)}
  \qquad\text{for every }f\in L^2(\rho).
\]
Let \(\gamma_{\mathrm{op}}(K)\) denote the supremum of \(\rho(K)\)
over all such measures. By \cite[Section~1.2, (8)--(10)]{CMPT},
\[
  c_1\gamma(K)\le\gamma_{\mathrm{op}}(K)\le C_1\gamma(K),
\]
where \(c_1,C_1>0\) are independent of \(K\). Choose an admissible measure
\(\nu\) on \(L\) with \(\nu(L)>\tfrac12\gamma_{\mathrm{op}}(L)\).
For every compact set \(J\subset\mathbb C\), the restriction
\(\nu|_J\) is admissible on \(J\). Hence
\[
  \nu(J)\le\gamma_{\mathrm{op}}(J)\le C_1\gamma(J).
\]
Hence, by inner regularity, for every Borel set \(F\subset\mathbb C\),
\begin{equation*}
  \tag{\ref*{cau:measure}a}\label{cau:borel-domination}
  \nu(F)
  =\sup_{\substack{J\subset F\\J\ \mathrm{compact}}}\nu(J)
  \le C_1\sup_{\substack{J\subset F\\J\ \mathrm{compact}}}\gamma(J)
  =C_1\gamma(F).
\end{equation*}

We extend the scalar weak-type estimate of Lemma \ref{cau:capacity} (1) to Hilbert-space-valued measures.

\begin{lemma}\label{cau:max}
There exists a constant \(C>0\), independent of the Hilbert space
\(V\), the measure \(\mu\), and the parameter \(t\), such that
\[
  \gamma\{\mathcal C_*\mu>t\}
  \le\frac{C}{t}\,|\mu|(\mathbb C), \
  \text{ for }t>0.
\]
\end{lemma}

\begin{proof}
Let \(v(\cdot)\) denote the Radon–Nikodým derivative of \(\mu\) with respect to its total variation measure \(|\mu|\), i.e., \(d\mu=v\,d|\mu|\). Then \(\|v\|_V=1\)
\(|\mu|\)-almost everywhere. Fix an orthonormal basis \(\{e_j\}\) of \(V\)
and let $P_n$ be the orthogonal projection defined as follows
\[
  P_nx=\sum_{j=1}^{\min\{n,\dim V\}}\langle x,e_j\rangle e_j,
  \text{ for }x\in V.
\]
If \(V\) is finite-dimensional, it suffices to take \(1\le n\le\dim V\)
in the Gaussian sums below, since \(P_n=I\) for \(n\ge\dim V\).
Write \(v_j(\cdot)=\langle v(\cdot),e_j\rangle\). For each \(n\), let \(g_1,\ldots,g_n\) be independent
standard complex Gaussians with \(\mathbb E|g_j|^2=1\), and \textcolor{black}{let}
\(d\xi_g=(\sum_{j=1}^n g_jv_j)\,d|\mu|\). Then
\[
  \mathbb E|\xi_g|(\mathbb C)
  =\int_{\mathbb C}\mathbb E\Bigl|\sum_{j=1}^n g_jv_j\Bigr|\,d|\mu|
  \le\int_{\mathbb C}\Bigl(\sum_{j=1}^n|v_j|^2\Bigr)^{1/2}d|\mu|
  \le|\mu|(\mathbb C).
\]

Fix \(t>0\) and \textcolor{black}{set} \(\widehat{E}_n=\{\mathcal C_*(P_n\mu)>t\}\).
For \(z\in \widehat{E}_n\), choose \(\varepsilon>0\) such that
\(s=\|P_n\mathcal C_\varepsilon\mu(z)\|_V>t\). The Gaussian identity
\[
  \mathcal C_\varepsilon\xi_g(z)
  =\sum_{j=1}^n g_j\bigl(\mathcal C_\varepsilon\mu(z)\bigr)_j
  \ \stackrel{\mathrm{law}}{=}\ s g_1
\]
implies
\[
  \mathbb P\{\mathcal C_*\xi_g(z)>t/2\}
  \ge\mathbb P\{|g_1|>t/(2s)\}
  =e^{-t^2/(4s^2)}\ge e^{-1/4}.
\]
By Lemma~\ref{cau:rational-truncations},
\[
  \mathcal C_*\xi_g(z)
  =\sup_{q\in\mathbb Q\cap(0,\infty)}
    \left|\sum_{j=1}^n g_j\bigl(\mathcal C_q\mu(z)\bigr)_j\right|
\]
is jointly measurable in \((z,g)\).

Let \(L\subset\{\mathcal C_*\mu>t\}\) be compact with \(\gamma(L)>0\).
By Lemma~\ref{cau:capacity}, choose a positive measure \(\nu\)
supported on \(L\) satisfying \eqref{cau:measure}.
For each fixed \(g\), \textcolor{black}{set}
\(E_g=\{z\in\mathbb C:\mathcal C_*\xi_g(z)>t/2\}\). Then
\[
  \nu(E_g)
  \overset{\text{\scriptsize\eqref{cau:borel-domination}}}{\le} C_1\gamma(E_g)
  \overset{\text{\scriptsize\eqref{cau:scalar}}}{\le}
  \frac{C}{t}\,|\xi_g|(\mathbb C).
\]
Taking expectations and applying Tonelli's theorem, we obtain
\[
\begin{aligned}
  e^{-1/4}\nu(\widehat{E}_n)
  &\le\int_{\mathbb C}\mathbb P\{\mathcal C_*\xi_g(z)>t/2\}\,d\nu(z)\\
  &=\mathbb E\,\nu\{\mathcal C_*\xi_g>t/2\}
   \le\frac Ct\,\mathbb E|\xi_g|(\mathbb C)
   \le\frac Ct\,|\mu|(\mathbb C).
\end{aligned}
\]
The bound is independent of \(n\). Since
\[
 \lim\limits_{n\to\infty} \|P_n\mathcal C_\varepsilon\mu(z)\|_V
  =\|\mathcal C_\varepsilon\mu(z)\|_V,
  \text{ for }\varepsilon>0,
\]
and \(P_n\le P_{n+1}\), we have
\[
  \widehat{E}_n\subset\widehat{E}_{n+1},\text{ and }
  \bigcup_{n\ge1}\widehat{E}_n=\{\mathcal C_*\mu>t\}.
\]
The support condition on \(\nu\) and continuity from below yield
\[
  c\gamma(L)\le\nu(L)
  =\lim\limits_{n\to\infty}\nu(\widehat{E}_n)
  \le\frac Ct\,|\mu|(\mathbb C).
\]
Taking the supremum over compact subsets, we obtain
\[
  \gamma\{\mathcal C_*\mu>t\}
  =\sup\left\{\gamma(L);\,L\subset\{\mathcal C_*\mu>t\} \text{ and }L\text{ is } \mathrm{compact}\right\}
  \le\frac Ct\,|\mu|(\mathbb C).
\]
\end{proof}

By combining Tolsa’s scalar principal value theorem \cite[Section 3, p. 1517]{TolsaICM} with the maximal estimation in Lemma 2.4, we obtain the following Hilbert-space-valued extension by finite-dimensional approximation.

\begin{lemma}\label{cau:pv}
The principal value
\[
  \mathcal C\mu(z)=\lim\limits_{\varepsilon\to0^+}\mathcal C_\varepsilon\mu(z)
\]
exists in \(V\) outside a Borel set of analytic capacity zero.
\end{lemma}

\begin{proof}
By Lemma~\ref{cau:rational-truncations}, the oscillation at zero is
\[
\begin{aligned}
  o_\mu(z)
  &:=\lim\limits_{k\to\infty}
     \sup_{0<\varepsilon,\delta<1/k}
     \|\mathcal C_\varepsilon\mu(z)-\mathcal C_\delta\mu(z)\|_V\\
  &=\lim\limits_{k\to\infty}
     \sup_{q,s\in\mathbb Q\cap(0,1/k)}
     \|\mathcal C_q\mu(z)-\mathcal C_s\mu(z)\|_V\\
  &=\inf_{k\ge1}
     \sup_{q,s\in\mathbb Q\cap(0,1/k)}
     \|\mathcal C_q\mu(z)-\mathcal C_s\mu(z)\|_V.
\end{aligned}
\]
The last expression is a countable infimum of countable suprema
of Borel functions. Thus \(o_\mu\) is a Borel function.

Set \(\mathcal A=\{o_\mu>0\}\).
Suppose that \(\gamma(\mathcal A)>0\).
By the definition of analytic capacity, there is a compact set
\(L\subset\mathcal A\) with \(\gamma(L)>0\).
By Lemma~\ref{cau:capacity}, choose a positive measure \(\nu\)
supported on \(L\) satisfying \eqref{cau:measure}. In particular,
\[
  0<c\gamma(L)\le\nu(L).
\]
It is therefore enough to prove that \(\nu(\mathcal A)=0\).

Let \(\{e_j\}\) and \(P_n\) be as in the proof of
Lemma~\ref{cau:max}, and set \(R_n=I-P_n\).
By the triangle inequality,
\[
  o_\mu\le o_{P_n\mu}+2\mathcal C_*(R_n\mu).
\]
At first, we show that \(o_{P_n\mu}=0\)
\(\nu\)-almost everywhere for every \(n\).
\textcolor{black}{Let} \(\mu_j=\langle\mu,e_j\rangle\) and \(Z_j=\{o_{\mu_j}>0\}\).
By Tolsa's principal value theorem \cite[Section~3, p.~1517]{TolsaICM}
and \eqref{cau:borel-domination},
\[
  \gamma(Z_j)=0,\text{ and } \nu(Z_j)\le C\gamma(Z_j)=0.
\]
Hence \(Z=\bigcup_j Z_j\) is \(\nu\)-null.
For \(z\notin Z\) and each \(n\), orthogonality gives
\begin{eqnarray*}
  &&\limsup\limits_{\varepsilon,\delta\to0^+}\|\mathcal C_\varepsilon(P_n\mu)(z)-\mathcal C_\delta(P_n\mu)(z)\|_V^2\\
  &&~~~=\limsup\limits_{\varepsilon,\delta\to0^+}\sum_{j=1}^{\min\{n,\dim V\}}
    |\mathcal C_\varepsilon\mu_j(z)-\mathcal C_\delta\mu_j(z)|^2
  =0.
\end{eqnarray*}
Thus \(o_{P_n\mu}(z)=0\) for every \(n\) and every \(z\notin Z\).

It remains to estimate \(\mathcal C_*(R_n\mu)\).
Write \(d\mu=v\,d|\mu|\), where \(\|v\|_V=1\)
almost everywhere with respect to \(|\mu|\).
Since \(o_{P_n\mu}=0\) \(\nu\)-almost everywhere,
by \eqref{cau:borel-domination} and Lemma~\ref{cau:max}, for every \(t>0\),
\[
\begin{aligned}
  \limsup\limits_{n\to\infty}\nu\{o_\mu>t\}
  &\le \limsup\limits_{n\to\infty}\nu\{\mathcal C_*(R_n\mu)>t/2\}\\
  & \le \limsup\limits_{n\to\infty}C\gamma\{\mathcal C_*(R_n\mu)>t/2\}\\
  &\le \limsup\limits_{n\to\infty}\frac{C}{t}\,|R_n\mu|(\mathbb C)
   = \limsup\limits_{n\to\infty}\frac{C}{t}\int_{\mathbb C}\|R_nv\|_V\,d|\mu|
  =0,
\end{aligned}
\]
where the last limit follows by Dominated Convergence Theorem, provided that
\(R_nv\to0\) and \(\|R_nv\|_V\le1\)
almost everywhere with respect to \(|\mu|\).
Consequently, \(\nu\{o_\mu>t\}=0\) for every \(t>0\).
Since
\[
  \mathcal A=\bigcup_{m\ge1}\{o_\mu>1/m\},
\]
by the continuity from below, one has
\[
  \nu(\mathcal A)
  =\lim\limits_{m\to\infty}\nu\{o_\mu>1/m\}=0.
\]
However \(L\subset\mathcal A\), it follows that
\[
  0<c\gamma(L)\le\nu(L)\le\nu(\mathcal A)=0,
\]
which is a contradiction. Hence \(\gamma(\mathcal A)=0\).

For \(z\notin\mathcal A\), we have \(o_\mu(z)=0\), and therefore
\[
  \lim\limits_{\varepsilon,\delta\to0^+}
  \|\mathcal C_\varepsilon\mu(z)-\mathcal C_\delta\mu(z)\|_V=0.
\]
Since \(V\) is complete, the Cauchy criterion gives the norm limit
\[
  \mathcal C\mu(z):=\lim\limits_{\varepsilon\to0^+}\mathcal C_\varepsilon\mu(z)\in V.
\]
\end{proof}
Using Lemma 2.4, we obtain a Hilbert-space-valued analogue of the scalar continuity estimate in  Akeroyd, Conway, and Yang \cite[Lemma 3.2]{AkeroydConwayYang}.
\begin{lemma}\label{cau:density}
Let \(a\in\mathbb C\) and \(a_0\in V\). Suppose that
\[
  \lim_{t\to0^+}\mathcal C_t\mu(a)=a_0 \quad\text{in }V,
  \qquad
  \lim_{t\to0^+}\frac{|\mu|(B(a,t))}{t}=0.
\]
For every \(\delta>0\) there exist \(r_0>0\) and Borel sets
\(E_r\subset B(a,r)\), \(0<r<r_0\), such that
\[
  \lim_{r\to0^+}\frac{\gamma(E_r)}{r}=0,
  \qquad
  \|\mathcal C\mu-a_0\|_{L^\infty(B(a,r)\setminus E_r,\,m_2;\,V)}
  \le\delta,
  \qquad 0<r<r_0.
\]
\end{lemma}

\begin{proof}
\textcolor{black}{Let} \(\theta(t)=|\mu|(B(a,t))/t\). For \(r>0\), write
\[
  \mu_0=\mu|_{\overline{B(a,2r)}},
  \qquad \mu_1=\mu-\mu_0.
\]
Then
\[
  |\mu_0|(\mathbb C)\le 3r\theta(3r)=o(r).
\]
By the definition of \(\mu_1\) and the hypothesis at \(a\),
\[
  \mathcal C\mu_1(a)=\mathcal C_{2r}\mu(a)
  \longrightarrow a_0 \qquad \text{as }r\downarrow0.
\]
If \(|z-a|<r\) and \(|w-a|>2r\), then \(|z-w|\ge |w-a|/2\).
It follows that
\[
  \left|\frac{1}{z-w}-\frac{1}{a-w}\right|
  =\frac{|z-a|}{|z-w|\,|a-w|}
  \le\frac{2r}{|w-a|^2}.
\]
For \(k\ge1\), set \(A_k=\{w:2^kr<|w-a|\le 2^{k+1}r\}\). Then
\[
\begin{aligned}
  \limsup_{r\to0^+}\sup_{z\in B(a,r)}
  \|\mathcal C\mu_1(z)-\mathcal C\mu_1(a)\|_V
  &\le \limsup_{r\to0^+}
       2r\sum\limits_{k\ge1}\frac{|\mu|(A_k)}{(2^kr)^2}\\
  &\le \limsup_{r\to0^+}
       8\sum\limits_{k\ge1}2^{-k}\theta(2^{k+2}r)
   =0.
\end{aligned}
\]
The last limit follows by dominated convergence, since \(\theta\)
is bounded on \((0,\infty)\) and \(\theta(t)\to0\) as \(t\downarrow0\).
Hence
\[
\begin{aligned}
  \sup_{z\in B(a,r)}\|\mathcal C\mu_1(z)-a_0\|_V
  &\le \sup_{z\in B(a,r)}
       \|\mathcal C\mu_1(z)-\mathcal C\mu_1(a)\|_V
       +\|\mathcal C_{2r}\mu(a)-a_0\|_V\\
  &\longrightarrow 0 \qquad \text{as }r\downarrow0.
\end{aligned}
\]

Fix \(\delta>0\) and \textcolor{black}{set}
\[
  E_r=B(a,r)\cap\{z:\mathcal C_*\mu_0(z)>\delta/2\}.
\]
By Lemma~\ref{cau:rational-truncations}, \(E_r\) is Borel, and
Lemma~\ref{cau:max} gives
\[
  \gamma(E_r)
  \le \frac{2C}{\delta}\,|\mu_0|(\mathbb C)
  =o(r).
\]
By Lemma~\ref{cau:area-lemma},
\(\mathcal C\mu=\mathcal C\mu_0+\mathcal C\mu_1\)
at \(m_2\)-almost every point of \(B(a,r)\). Thus, for
\(m_2\)-almost every \(z\in B(a,r)\setminus E_r\),
\[
  \|\mathcal C\mu(z)-a_0\|_V
  \le \mathcal C_*\mu_0(z)
     +\sup_{w\in B(a,r)}\|\mathcal C\mu_1(w)-a_0\|_V
  \le \frac{\delta}{2}+o(1).
\]
The last bound is at most \(\delta\) for all sufficiently small \(r>0\).
\end{proof}
The following lemma is a re-scaled dual formulation of the point evaluation estimate of Aleman, Richter and Sundberg\cite[p. 458, Lemma B]{ARS}.

\begin{lemma}\label{cau:representing-density}
There exist constants \(c_0,C_0>0\), independent on the choice of the center \(a\),
the radius \(r\), and the Borel set \(E\), with the following property.
For \(E\subset B(a,r)\) and \(\gamma(E)<c_0r\), set
\(G=B(a,r)\setminus E\). Then there is an \(h\in L^\infty(G,m_2)\)
such that for $p\in\mathbb C[z],$
\begin{equation}\label{cau:represent}
  \int_G ph\,dm_2=p(a),\,
  \|h\|_\infty\le C_0r^{-2}, \text{ and }
  \|h\|_1\le\pi C_0.
\end{equation}
\end{lemma}

\begin{proof}
Since \(\gamma((E-a)/r)=r^{-1}\gamma(E)<c_0\),
applying \cite[p.~458, Lemma~B]{ARS} to the map
\(\zeta\mapsto p(a+r\zeta)\) implies
\[
  |p(a)|
  \le C_0\int_{B(0,1)\setminus((E-a)/r)}
       |p(a+r\zeta)|\,dm_2(\zeta)
  =C_0r^{-2}\int_G|p(z)|\,dm_2(z).
\]
Thus \(p|_G\mapsto p(a)\) is a well-defined bounded linear functional
on the polynomial restrictions in \(L^1(G,m_2)\), with norm at most
\(C_0r^{-2}\).
By Hahn--Banach Theorem and \((L^1)^*=L^\infty\), there is an
\(h\in L^\infty(G,m_2)\) such that for any $p\in\mathbb C[z]$
\[
  \int_G ph\,dm_2=p(a), \text{ and }
  \|h\|_\infty\le C_0r^{-2}.
\]
Consequently,
\[
  \|h\|_1
  \le m_2(G)\cdot\|h\|_\infty
  \le \pi r^2 C_0r^{-2}
  =\pi C_0.
\]
\end{proof}

For a finite positive Borel measure \(\mu\) on \(\mathbb C\), its Riesz \(1\)-energy is defined by $$I_1(\mu):=\int_{\mathbb C} I_\mu(z)\,d\mu(z),$$ where \(I_\mu\) is the potential introduced in Lemma \ref{cau:area-lemma}.
The following lemma uses Egorov’s theorem to restrict a measure of finite Riesz \(1\)-energy to a compact set of positive measure on which the restriction has uniformly small local potentials.

\begin{lemma}\label{cau:energy}
Let \(\mu\) be a probability measure supported on a compact set
\(L\subset\mathbb C\), and suppose that
\[
  I_1(\mu)
  =\iint_{L\times L}\frac{d\mu(z)\,d\mu(\zeta)}{|z-\zeta|}
  <\infty.
\]
There is a compact set \(F\subset L\) with \(\mu(F)>0\) such that,
for \(\tau=\mu|_F\),
\[
  b_\tau(r):=\sup_{w\in\mathbb C}
  \int_{B(w,r)}\frac{d\tau(z)}{|z-w|}
  \rightarrow 0, \text{ as }r\to0^+.
\]
In particular,
\[
  C_\tau:=\sup_{w\in\mathbb C}\int_F\frac{d\tau(z)}{|z-w|}
  \le \inf_{r>0}\left(b_\tau(r)+\frac{\tau(F)}{r}\right)<\infty,
\]
and \(\tau(J)\le C_\tau\,\gamma(J)\) for any compact subset \(J\subset L\).
\end{lemma}

\begin{proof}
Since \(\int_L I_\mu\,d\mu=I_1(\mu)<\infty\), we have
$$\mu\{a\in L:I_\mu(a)=\infty\}=0.$$
If \(a\) were an atom of \(\mu\), then \(I_\mu(a)=\infty\), so
\[
  I_1(\mu)\ge\int_{\{a\}} I_\mu\,d\mu=\infty,
\]
which is a contradiction. Thus \(\mu\) has no atoms.
For \(r>0\), set
\[
  u_r(a)=\int_{B(a,r)}\frac{d\mu(z)}{|z-a|}.
\]
At each \(a\in L\) with \(I_\mu(a)<\infty\), by Dominated
Convergence Theorem,
\[
 \limsup\limits_{n\to\infty} u_{2^{-n}}(a)
  =\limsup\limits_{n\to\infty}\int_{0<|z-a|<2^{-n}}\frac{d\mu(z)}{|z-a|}
  =0.
\]
By Egorov's theorem, choose a Borel set \(F_0\subset L\) of positive
\(\mu\)-measure on which \(u_{2^{-n}}\to0\) uniformly.
Inner regularity then gives a compact \(F\subset F_0\) with
\(\mu(F)>0\). Hence
\[
  0\le\limsup\limits_{n\to\infty}\sup_{a\in F}u_{2^{-n}}(a)
  \le\limsup\limits_{n\to\infty}\sup_{a\in F_0}u_{2^{-n}}(a)
  =0.
\]
The monotonicity of \(u_r(a)\) in \(r\) extends this convergence
to arbitrary radii:
\[
  0\le \limsup\limits_{r\to0^+}\sup_{a\in F}u_r(a)
  \le \inf_{n\ge1}\sup_{a\in F}u_{2^{-n}}(a)=0.
\]
Let \(\tau=\mu|_F\).
We next estimate its potential at arbitrary centers.
Fix \(w\in\mathbb C\). By compactness of \(F\), choose \(a\in F\)
with \(|w-a|=\operatorname{dist}(w,F)\). Then
\[
  |z-a|\le |z-w|+|w-a|\le2|z-w|
  ,\text{ for } z\in F.
\]
Therefore, \(F\cap B(w,r)\subset F\cap B(a,2r)\). Consequently,
\[
\begin{aligned}
  \int_{B(w,r)}\frac{d\tau(z)}{|z-w|}
  &=\int_{F\cap B(w,r)}\frac{d\mu(z)}{|z-w|}\\
  &\le2\int_{F\cap B(a,2r)}\frac{d\mu(z)}{|z-a|}\\
  &\le2u_{2r}(a)
   \le2\sup_{\zeta\in F}u_{2r}(\zeta).
\end{aligned}
\]
Taking the supremum over \(w\), we obtain
\[
  \limsup\limits_{r\to0^+}b_\tau(r)\le \limsup\limits_{r\to0^+}\sup_{a\in F}2u_{2r}(a)
  =0.
\]
%We may therefore choose \(r_0>0\) with \(b_\tau(r_0)<\infty\).
For every \(r>0\),
\[
\begin{aligned}
  C_\tau
  =\sup_{w\in\mathbb C}
    \left(
      \int_{F\cap B(w,r)}\frac{d\tau(z)}{|z-w|}
      +\int_{F\setminus B(w,r)}\frac{d\tau(z)}{|z-w|}
    \right)
  \le b_\tau(r)+\frac{\tau(F)}{r}.
\end{aligned}
\]
Since \(b_\tau(r)\to0\) as \(r\to0^+\), it follows that
\[
  C_\tau\le\inf_{r>0}
  \left(b_\tau(r)+\frac{\tau(F)}{r}\right)<\infty.
\]
Since \(\tau(F)=\mu(F)>0\), we also have \(C_\tau>0\).

Finally, for compact \(J\subset L\), set
\[
  f_J(z)=\frac{1}{C_\tau}\int_J\frac{d\tau(w)}{z-w},
  \,\, z\in\mathbb C\setminus J.
\]
Then \(f_J\) is holomorphic off \(J\) and  \(f_J(\infty)=0\).
By the definition of \(C_\tau\), \(\|f_J\|_\infty\le1\).
Thus \(f_J\) is admissible for \(\gamma(J)\), and
\[
  \gamma(J)\ge |f_J'(\infty)|
  =\lim\limits_{z\to\infty}|z f_J(z)|
  =\frac{\tau(J)}{C_\tau}.
\]
\end{proof}

The following consequence of Lemmas~\ref{cau:density} and
\ref{cau:representing-density} will allow the representing density to be
supported away from its evaluation point.

\begin{lemma}\label{cau:annular-density}
Let $U$ be an open disk and let $a\in U$ satisfy $I_\mu(a)<\infty$.
For every $\delta>0$ there exist $0<\rho<r$ with
$\overline{B(a,r)}\subset U$, a Borel set
\[
  G\subset B(a,r)\setminus\overline{B(a,\rho)},
\]
and a function $h\in L^\infty(G,m_2)$ such that
\[
  \|\mathcal C\mu(z)-\mathcal C\mu(a)\|_V\le\delta
  \quad\text{for $m_2$-almost every }z\in G,
\]
and for any $p\in\mathbb C[z]$
\[
  p(a)=\int_G p(z)h(z)\,dm_2(z).
\]
\end{lemma}
\begin{proof}
Since $I_\mu(a)<\infty$, we have $|\mu|(\{a\})=0$, and the singular integral
$\int_{\mathbb C}(a-w)^{-1}\,d\mu(w)$ converges to some $y$ absolutely.
By Dominated Convergence Theorem,
\[
  \max\left\{
    \|\mathcal C_t\mu(a)-y\|_V,\,
    \frac{|\mu|(B(a,t))}{t}
  \right\}
  \le \int_{0<|w-a|\le t}\frac{d|\mu|(w)}{|a-w|}
  \longrightarrow 0,
  \text{ as } t\to0^+.
\]
In particular, $\mathcal C\mu(a)=y$.

Fix \(\delta>0\). By Lemma~\ref{cau:density} with \(a_0=y\),
for all sufficiently small \(r>0\) we can choose Borel sets
\(\widehat E_r\subset B(a,r)\) such that
\[
  \gamma(\widehat E_r)=o(r)\qquad\text{as }r\rightarrow0.
\]
Moreover,
\[
  \|\mathcal C\mu(z)-\mathcal C\mu(a)\|_V\le\delta
  \qquad\text{for $m_2$-almost every }z\in B(a,r)\setminus \widehat{E}_r.
\]
For $0<\theta<1$, let
\[
  \widehat{E'}_r=\widehat{E}_r\cup\overline{B(a,\theta r)}\subset B(a,r).
\]
Let $L\subset \widehat{E'}_r$ be compact with $\gamma(L)>0$, and let $\nu$
be the measure given by Lemma~\ref{cau:capacity}\textup{(2)} for $L$.
For every $\varepsilon>0$ we have
$L\subset \widehat{E}_r\cup B(a,\theta r+\varepsilon)$. Thus, with $c,C$
as in \eqref{cau:measure} and using \eqref{cau:borel-domination}
for \(\widehat{E}_r\),
\[
  c\gamma(L)
  \le \nu(L)
  \le \nu(\widehat{E}_r)+\nu\bigl(B(a,\theta r+\varepsilon)\bigr)
  \le C\gamma(\widehat{E}_r)+\theta r+\varepsilon.
\]
Letting $\varepsilon\to0^+$ gives
$c\gamma(L)\le C\gamma(\widehat{E}_r)+\theta r$. Taking the supremum over
compact $L\subset \widehat{E'}_r$, we obtain
\[
  c\gamma(\widehat{E'}_r)\le C\gamma(\widehat{E}_r)+\theta r.
\]
Let $c_0$ be the constant in Lemma~\ref{cau:representing-density}.
Choose $0<\theta<1$ satisfying $2\theta<cc_0$. Since $\gamma(\widehat{E}_r)=o(r)$,
\[
  \limsup_{r\to0^+}\frac{\gamma(\widehat{E'}_r)}{r}
  \le \frac{\theta}{c}<\frac{c_0}{2}.
\]
Choose $r>0$ sufficiently small that $\overline{B(a,r)}\subset U$
and $\gamma(\widehat{E'}_r)<c_0r$. \textcolor{black}{Set} $\rho=\theta r$ and
$G=B(a,r)\setminus \widehat{E'}_r$. Then
\[
  G\subset B(a,r)\setminus\overline{B(a,\rho)},
  \qquad \operatorname{dist}(a,G)\ge\rho>0.
\]
Since $G\subset B(a,r)\setminus \widehat{E}_r$, the estimate on
$B(a,r)\setminus \widehat{E}_r$ implies that
\[
  \|\mathcal C\mu(z)-\mathcal C\mu(a)\|_V\le\delta
  \qquad\text{for $m_2$-almost every }z\in G.
\]
Finally, Lemma~\ref{cau:representing-density}, applied with
\(E=\widehat{E'}_r\), gives \(h\in L^\infty(G,m_2)\) such that
\[
  \int_G p(z)h(z)\,dm_2(z)=p(a),
  \qquad p\in\mathbb C[z].
\]
\end{proof}

\section{Spectral integrals and local polynomial bounds}\label{sec:model}

Let \(H\) be a separable complex Hilbert space, let
\(N\in\mathcal B(H)\) be normal, and let \(K\in S_1(H)\).
\textcolor{black}{Set} \(T=N+K\), \(\Lambda=\sigma(N)\), and
\(\mathcal E=H\oplus H\), and let
\(E\) be the spectral measure of \(N\).

\begin{lemma}\label{model:factor}
There are Hilbert--Schmidt operators
\(A\colon\mathcal E\to H\) and \(B\colon H\to\mathcal E\) with
the following properties.
\begin{enumerate}
\item
\(AB=K\) and \(\overline{A\mathcal E}=H\).

\item
The formula \(\nu_A(F)=\|E(F)A\|_2^2\), for Borel
\(F\subset\mathbb C\), defines a finite positive measure with
\[
  \nu_A(\mathbb C)=\|A\|_2^2,\qquad
  \operatorname{supp}\nu_A=\Lambda,
\]
where \(\|\cdot\|_2\) denotes the Hilbert--Schmidt norm.

\item
The formula \(\mu(F)=BE(F)A\), for Borel \(F\subset\mathbb C\),
defines a countably additive \(S_2(\mathcal E)\)-valued measure with
\[
  |\mu|(\mathbb C)\le\|B\|_2\|A\|_2.
\]
For every bounded Borel function \(\varphi\) on \(\Lambda\),
\begin{equation}\label{model:borel-calculus}
  B\varphi(N)A=\int_\Lambda\varphi\,d\mu.
\end{equation}
\end{enumerate}
\end{lemma}

\begin{proof}
\textup{(1)}
Write $K=V_0|K|$ and \textcolor{black}{set} $S_0=|K|^{1/2}$. Choose $J_0\in S_2(H)$
with dense range and set
\[
  A(u,v)=V_0S_0u+J_0v,\qquad Bx=(S_0x,0).
\]
Since $\|S_0\|_2^2=\operatorname{tr}|K|<\infty$, both $A$ and $B$
are Hilbert--Schmidt. Also $AB=V_0S_0^2=K$ and
$J_0H\subset A\mathcal E$.

\smallskip
\noindent\textup{(2)}
For an orthonormal basis $\textcolor{black}{\{e_n\}}$ of $\mathcal E$,
\[
  \nu_A(F)=\sum_n\langle E(F)Ae_n,Ae_n\rangle .
\]
Hence $\nu_A$ is a positive measure of mass $\|A\|_2^2$. For open $U$,
\[
  \nu_A(U)=0
  \ \Longleftrightarrow\ E(U)A=0
  \ \Longleftrightarrow\ E(U)=0
  \ \Longleftrightarrow\ U\cap\Lambda=\varnothing .
\]
Thus $\operatorname{supp}\nu_A=\Lambda$.

\smallskip
\noindent\textup{(3)}
Let $(F_j)$ be pairwise disjoint Borel sets and \textcolor{black}{let} $E_j=E(F_j)$.
By orthogonality,
\[
\begin{aligned}
  \sum_j\|BE_jA\|_2
  &\le \sum_j\|E_jB^*\|_2\,\|E_jA\|_2\\
  &\le \left(\sum_j\|E_jB^*\|_2^2\right)^{1/2}
      \left(\sum_j\|E_jA\|_2^2\right)^{1/2}
  \le \|B\|_2\|A\|_2 .
\end{aligned}
\]
If $F=\bigcup_jF_j$ and $R_m=\bigcup_{j>m}F_j$, then
\[
  \left\|\mu(F)-\sum_{j=1}^m\mu(F_j)\right\|_2
  =\|BE(R_m)A\|_2
  \le\|B\|_{\mathrm{op}}\,\nu_A(R_m)^{1/2}
  \xrightarrow{m\to\infty}0 .
\]
Thus \(\mu\) is countably additive. Let \(\mathcal P\) range over
the finite Borel partitions of \(\mathbb C\). Then
\[
  |\mu|(\mathbb C)
  =\sup_{\mathcal P}\sum_{F\in\mathcal P}\|BE(F)A\|_2
  \le\|B\|_2\|A\|_2.
\]
For a bounded Borel function \(\varphi\) on \(\Lambda\), choose
simple functions satisfying
\[
  \varphi_n=\sum_j c_{nj}\mathbf1_{F_{nj}},\qquad
  \|\varphi_n-\varphi\|_\infty\xrightarrow{n\to\infty}0.
\]
Then
\[
  B\varphi_n(N)A
  =\sum_j c_{nj}BE(F_{nj})A
  =\sum_j c_{nj}\mu(F_{nj})
  =\int_\Lambda\varphi_n\,d\mu.
\]
Consequently,
\[
\begin{aligned}
  \left\|B\varphi(N)A-\int_\Lambda\varphi\,d\mu\right\|_2
  &\le
  \bigl(\|B\|_{\mathrm{op}}\|A\|_2+|\mu|(\mathbb C)\bigr)
  \|\varphi-\varphi_n\|_\infty\\
  &\xrightarrow{n\to\infty}0,
\end{aligned}
\]
which proves \eqref{model:borel-calculus}.
\end{proof}

Fix $A,B,\nu_A,\mu$ as in Lemma~\ref{model:factor}, so that $T=N+AB$. Write
$d\mu=M_0\,d|\mu|$, with $M_0$ taking values in $S_2(\mathcal E)$,
and identify $\mathcal B(\mathbb C,\mathcal E)$ with $\mathcal E$.

\begin{lemma}\label{model:spectral-integral}
Let \(Y\in\{\mathbb C,\mathcal E\}\). For a finite Borel partition
\((F_j)\) of \(\Lambda\) and \(g_j\in\mathcal B(Y,\mathcal E)\), define
\[
  \Gamma_A\left(\sum_j\mathbf1_{F_j}g_j\right)
  =\sum_j E(F_j)Ag_j.
\]
Then the following hold.
\begin{enumerate}
\item
The rule extends by uniform approximation to a bounded linear map
\[
  \Gamma_A\colon C(\Lambda,\mathcal B(Y,\mathcal E))
  \longrightarrow S_2(Y,H).
\]
For every Borel set \(F\subset\mathbb C\),
\begin{equation}\label{model:spectral-bound}
  \|E(F)\Gamma_A(g)\|_2^2
  \le\int_{F\cap\Lambda}\|g(w)\|_{\mathrm{op}}^2\,d\nu_A(w).
\end{equation}

\item
For \(g\in C(\Lambda,\mathcal B(\mathcal E))\),
\(b\in C(\Lambda,\mathcal E)\), and \(v\in\mathcal E\),
\begin{equation}\label{model:spectral-identities}
\begin{gathered}
  \Gamma_A(gv)=\Gamma_A(g)v,\qquad
  N\Gamma_A(b)=\Gamma_A(wb),\\[3pt]
  B\Gamma_A(b)=\int_\Lambda d\mu(w)\,b(w)
  =\int_\Lambda M_0(w)b(w)\,d|\mu|(w).
\end{gathered}
\end{equation}

\item
For constant fields,
\(\Gamma_A(I_{\mathcal E})=A\) and
\(\Gamma_A(v)=Av\) for \(v\in\mathcal E\).
\end{enumerate}
\end{lemma}
\begin{proof}
We use the spectral calculus in
\cite[Theorems~12.21 and~12.23]{RudinFA}.\\
\textup{(1)}
For simple \(g\), orthogonality gives
\[
\begin{aligned}
  \|E(F)\Gamma_A(g)\|_2^2
  &=\sum_j\|E(F\cap F_j)Ag_j\|_2^2\\
  &\le\sum_j\|g_j\|_{\mathrm{op}}^2\,\nu_A(F\cap F_j)
   =\int_{F\cap\Lambda}\|g(w)\|_{\mathrm{op}}^2\,d\nu_A(w).
\end{aligned}
\]
Taking \(F=\Lambda\) gives
\(\|\Gamma_A(g)\|_2\le\|A\|_2\|g\|_\infty\).
The extension and \eqref{model:spectral-bound} follow by uniform
approximation.

\smallskip
\noindent\textup{(2)}
Let \(g,b,s\) be constant on each member of a finite Borel
partition \((F_j)\) of \(\Lambda\), with respective values
\(g_j\in\mathcal B(\mathcal E)\), \(b_j\in\mathcal E\), and
\(s_j\in\mathbb C\) on \(F_j\). Then
\[
  \Gamma_A(gv)
  =\sum_j E(F_j)Ag_jv
  =\Gamma_A(g)v,
\]
\[
  s(N)\Gamma_A(b)
  =\sum_{i,j}s_iE(F_i)E(F_j)Ab_j
  =\sum_j s_jE(F_j)Ab_j
  =\Gamma_A(sb),
\]
and
\[
  B\Gamma_A(b)
  =\sum_j BE(F_j)Ab_j
  =\sum_j\mu(F_j)b_j
  =\int_\Lambda M_0(w)b(w)\,d|\mu|(w).
\]

Now let \(g,b\) be continuous. Choose simple fields \(g_n,b_n\)
converging uniformly to \(g,b\), respectively. By \textup{(1)},
\[
  \Gamma_A(gv)
  =\lim_{n\to\infty}\Gamma_A(g_nv)
  =\lim_{n\to\infty}\Gamma_A(g_n)v
  =\Gamma_A(g)v.
\]
Choose simple scalar functions \(s_n\) on \(\Lambda\) such that
\[
  \sup_{w\in\Lambda}|s_n(w)-w|\xrightarrow{n\to\infty}0.
\] Then \(s_nb_n\to wb\) uniformly, and
\[
  \|s_n(N)-N\|_{\mathrm{op}}
  \le\sup_{w\in\Lambda}|s_n(w)-w|
  \xrightarrow{n\to\infty}0.
\]
Hence
\[
  N\Gamma_A(b)
  =\lim_{n\to\infty}s_n(N)\Gamma_A(b_n)
  =\lim_{n\to\infty}\Gamma_A(s_nb_n)
  =\Gamma_A(wb).
\]
Finally, since
\[
  \left\|\int_\Lambda M_0(w)(b_n(w)-b(w))\,d|\mu|(w)\right\|
  \le |\mu|(\Lambda)\|b_n-b\|_\infty
  \xrightarrow{n\to\infty}0,
\]
we obtain
\[
  B\Gamma_A(b)
  =\lim_{n\to\infty}B\Gamma_A(b_n)
  =\int_\Lambda M_0(w)b(w)\,d|\mu|(w).
\]
This proves \eqref{model:spectral-identities}.

\smallskip
\noindent\textup{(3)}
This follows from \(E(\Lambda)=I_H\).
\end{proof}

For $r>0$ \textcolor{black}{set}
\[
  k_r(u)=\frac{\overline u}{\max\{r,|u|\}^2},\qquad
  \widetilde{\mathcal C}_r\mu(z)=\int_\Lambda k_r(z-w)\,d\mu(w).
\]

\begin{lemma}\label{model:regularized-transform}
Define
\[
  D(z)=
  \begin{cases}
    \displaystyle\lim\limits_{\varepsilon\to0^+}\mathcal C_\varepsilon\mu(z),
      &\text{if the limit exists in }S_2(\mathcal E),\\
    0,&\text{otherwise}.
  \end{cases}
\]
Then the following hold.
\begin{enumerate}
\item
The field \(D\colon\mathbb C\to S_2(\mathcal E)\) is Borel, and
\(\Delta(z)=I_{\mathcal E}-D(z)\) is Fredholm with
\(\operatorname{ind}\Delta(z)=0\) for every \(z\in\mathbb C\).

\item
For \(r>0\),
\begin{equation}\label{model:regularization}
  \widetilde{\mathcal C}_r\mu(z)
  =\frac{2}{r^2}\int_0^r t\,\mathcal C_t\mu(z)\,dt,
  \qquad
  0\le1-uk_r(u)\le\mathbf1_{\{|u|<r\}}.
\end{equation}

\item
At every point \(z\) of principal value convergence,
\begin{equation}\label{model:regularized-limit}
  \lim\limits_{r\to0^+}
  \|\widetilde{\mathcal C}_r\mu(z)-D(z)\|_2=0.
\end{equation}
\end{enumerate}
\end{lemma}

\begin{proof}
\textup{(1)}
The Borel measurability of \(D\) follows from the proof of
Lemma~\ref{cau:pv} and the definition of \(D\). Since
\(D(z)\in S_2(\mathcal E)\subset\mathcal K(\mathcal E)\),
\cite[Theorems~4.23 and~4.25(a)]{RudinFA} gives that
\(\Delta(z)=I_{\mathcal E}-D(z)\) is Fredholm with
\(\operatorname{ind}\Delta(z)=0\).

\smallskip
\noindent\textup{(2)}
The disk Cauchy transform formula
\cite[Section~4.3.2, (4.38)]{AIM}, after scaling and normalization,
gives
\[
  k_r(u)=\frac{1}{\pi r^2}
         \int_{B(0,r)}\frac{dm_2(v)}{u-v}.
\]
\textcolor{black}{Set} \(q_t(u)=1/u\) for \(|u|>t\) and \(q_t(u)=0\) otherwise.
Direct integration gives
\[
  \frac{2}{r^2}\int_0^r t q_t(u)\,dt=k_r(u),\qquad
  \frac{2}{r^2}\int_0^r t|q_t(u)|\,dt
  =|k_r(u)|\le r^{-1}.
\]
Since \(|\mu|(\Lambda)<\infty\), Fubini's theorem yields
\[
  \widetilde{\mathcal C}_r\mu(z)
  =\int_\Lambda k_r(z-w)\,d\mu(w)
  =\frac{2}{r^2}\int_0^r t\,\mathcal C_t\mu(z)\,dt.
\]
Moreover,
\[
  1-uk_r(u)
  =\left(1-\frac{|u|^2}{r^2}\right)\mathbf1_{\{|u|<r\}}
  \in\bigl[0,\mathbf1_{\{|u|<r\}}\bigr].
\]

\smallskip
\noindent\textup{(3)}
At each point of principal value convergence, \textup{(2)} gives
\[
  \|\widetilde{\mathcal C}_r\mu(z)-D(z)\|_2
  \le\sup_{0<t<r}\|\mathcal C_t\mu(z)-D(z)\|_2
  \xrightarrow{r\to0^+}0.
\]
\end{proof}

\begin{lemma}\label{model:averaging}
Let \(U\subset\mathbb C\) be a bounded open disk. For a bounded
Bochner-measurable field \(b\) on \(U\), with values in
\(\mathcal E\) or \(\mathcal B(\mathcal E)\), \textcolor{black}{set}
\[
  G_b(w)=\int_U\frac{b(z)}{z-w}\,dm_2(z).
\]
Then the following hold.
\begin{enumerate}
\item
The field \(G_b\) is norm-continuous on \(\mathbb C\), and
\[
  \|G_b\|_\infty\le c_U\|b\|_\infty.
\]

\item
For \(\mathcal E\)-valued \(b\), \textcolor{black}{set}
\(Jb=\Gamma_A(G_b|_\Lambda)\). Then
\begin{equation}\label{model:averaging-bound}
  \|Jb\|
  \le\left(\int_\Lambda\|G_b(w)\|^2\,d\nu_A(w)\right)^{1/2}
  \le\|A\|_2\|G_b\|_\infty
  \le c_U\|A\|_2\|b\|_\infty,
\end{equation}
and
\[
  BJb=\int_U D(z)b(z)\,dm_2(z).
\]
Moreover,
\begin{equation}\label{model:intertwine}
  TJb=J(zb)-A\int_U\Delta(z)b(z)\,dm_2(z).
\end{equation}
\end{enumerate}
\end{lemma}

\begin{proof}
\textup{(1)}
For \(r>0\), \textcolor{black}{set}
\[
  G_{b,r}(w)=\int_U k_r(z-w)b(z)\,dm_2(z).
\]
Since \(k_r\) is bounded and continuous, dominated convergence
shows that \(G_{b,r}\) is norm-continuous. Also,
\[
  \|G_b-G_{b,r}\|_\infty
  \le\|b\|_\infty
      \int_{|u|<r}
      \left(\frac1{|u|}-\frac{|u|}{r^2}\right)dm_2(u)
  \le2\pi r\,\|b\|_\infty
  \xrightarrow{r\to0^+}0.
\]
Thus \(G_b\) is norm-continuous. Its bound follows from
\eqref{cau:area}.

\smallskip
\noindent\textup{(2)}
The estimate \eqref{model:averaging-bound} follows from
\eqref{model:spectral-bound} with \(F=\Lambda\). Furthermore,
\[
  \int_\Lambda\int_U
  \frac{\|b(z)\|}{|z-w|}\,dm_2(z)\,d|\mu|(w)
  \le c_U\|b\|_\infty|\mu|(\mathbb C)<\infty.
\]
By Lemma~\ref{cau:area-lemma}, the absolutely convergent Cauchy
integral equals \(D(z)\) for \(m_2\)-almost every \(z\).
Hence \eqref{model:spectral-identities} and Fubini's theorem give
\[
\begin{aligned}
  BJb
  &=\int_\Lambda d\mu(w)\,G_b(w)\\
  &=\int_U\left(\int_\Lambda\frac{d\mu(w)}{z-w}\right)b(z)\,dm_2(z)
   =\int_U D(z)b(z)\,dm_2(z).
\end{aligned}
\]
Since
\[
  G_{zb}(w)-wG_b(w)=\int_U b(z)\,dm_2(z),
\]
Lemma~\ref{model:spectral-integral} yields
\[
  J(zb)-NJb
  =\Gamma_A(G_{zb}-wG_b)
  =A\int_U b(z)\,dm_2(z).
\]
Combining these identities with \(T=N+AB\) and \(\Delta=I-D\),
we obtain
\[
\begin{aligned}
  TJb
  &=NJb+ABJb\\
  &=J(zb)-A\int_U b(z)\,dm_2(z)
          +A\int_U D(z)b(z)\,dm_2(z)\\
  &=J(zb)-A\int_U\Delta(z)b(z)\,dm_2(z).
\end{aligned}
\]
\end{proof}

\begin{lemma}\label{model:regular}
Let \(U\subset\mathbb C\) be a bounded open disk and let
\(a\in U\cap\Lambda\). Suppose that
\[
  \lim\limits_{r\to0^+}\|\mathcal C_r\mu(a)-D(a)\|_2=0,
  \qquad
  \lim\limits_{r\to0^+}\frac{|\mu|(B(a,r))}{r}=0.
\]
If \(L=\Delta(a)=I-D(a)\) is invertible on \(\mathcal E\), then
there exist \(x\in H\setminus\{0\}\) and \(C_x<\infty\) such that
\[
  \|p(T)x\|\le C_x\sup_{z\in U}|p(z)|
  \qquad(p\in\mathbb C[z]).
\]
\end{lemma}

\begin{proof}
Let \(c_0,C_0\) be as in Lemma~\ref{cau:representing-density}.
\textcolor{black}{Set} \(\rho=(1+2\pi C_0)^{-1}\) and
\(\delta=\rho/\|L^{-1}\|_{\mathrm{op}}\).
By Lemma~\ref{cau:density}, we may choose \(r>0\) and a Borel set
\(E_r\subset B(a,r)\) such that \(\overline{B(a,r)}\subset U\),
\(\gamma(E_r)<c_0r\), and
\[
  \operatorname*{ess\,sup}_{z\in G}\|\Delta(z)-L\|_2
  =\operatorname*{ess\,sup}_{z\in G}\|D(z)-D(a)\|_2
  \le\delta,
  \qquad G=B(a,r)\setminus E_r.
\]
\textcolor{black}{Let} \(W(z)=L^{-1}(\Delta(z)-L)\). Almost everywhere on \(G\),
\[
  \|W(z)\|_{\mathrm{op}}
  \le\|L^{-1}\|_{\mathrm{op}}\delta=\rho<1.
\]
Thus \(\Delta(z)=L(I+W(z))\) is invertible almost everywhere on \(G\).
The Neumann series gives
\[
  Q(z):=\Delta(z)^{-1}L
  =(I+W(z))^{-1}
  =\sum_{n=0}^{\infty}(-W(z))^n,
\]
and
\[
  \|Q(z)-I\|_{\mathrm{op}}
  \le\sum_{n=1}^{\infty}\rho^n
  =\frac{\rho}{1-\rho}
  =\frac{1}{2\pi C_0}.
\]
Set \(Q=I\) on the exceptional null set.

By Lemma~\ref{cau:representing-density}, choose
\(h\in L^\infty(G,m_2)\) such that
\[
  \int_G ph\,dm_2=p(a),\quad p\in\mathbb C[z],\qquad
  \|h\|_1\le\pi C_0.
\]
Define
\[
  \Phi(z)=
  \begin{cases}
    (z-a)h(z)Q(z),&z\in G,\\
    0,&z\in U\setminus G.
  \end{cases}
\]
Since \(D\) is Borel in the Hilbert--Schmidt norm, the Neumann series
shows that \(Q\) is measurable in the operator norm. For almost every
\(z\in G\),
\[
  W(z)=L^{-1}\bigl(D(a)-D(z)\bigr)\in S_2(\mathcal E),
  \qquad
  Q(z)-I=-Q(z)W(z)\in S_2(\mathcal E).
\]
Thus \(\Phi(z)=(z-a)h(z)Q(z)\) belongs to
\(\mathbb C I+S_2(\mathcal E)\) almost everywhere on \(G\).
This space is separable in the operator norm, and \(\Phi=0\) off \(G\).
Hence \(\Phi\) is Bochner-measurable on \(U\). Moreover,
\[
  \|\Phi\|_\infty
  \le r\|h\|_\infty\|Q\|_{L^\infty(G)}
  \le\frac{r\|h\|_\infty}{1-\rho}<\infty.
\]
Hence Lemma~\ref{model:averaging}\textup{(1)} applies to \(\Phi\), and
\(G_\Phi\) is norm-continuous.

Taking \(p=1\) in \eqref{cau:represent} gives
\(\int_G h\,dm_2=1\). Consequently,
\[
  \|G_\Phi(a)-I\|_{\mathrm{op}}
  =\left\|\int_G h(z)(Q(z)-I)\,dm_2(z)\right\|_{\mathrm{op}}
  \le\frac{\|h\|_1}{2\pi C_0}\le\frac12.
\]
By norm continuity, choose an open disk \(U_a\ni a\) such that
\(\sup_{w\in U_a}\|G_\Phi(w)-G_\Phi(a)\|_{\mathrm{op}}\le1/4\).
Then
\[
  \sup_{w\in U_a}\|G_\Phi(w)-I\|_{\mathrm{op}}
  \le\sup_{w\in U_a}\|G_\Phi(w)-G_\Phi(a)\|_{\mathrm{op}}
     +\|G_\Phi(a)-I\|_{\mathrm{op}}
  \le\frac14+\frac12=\frac34.
\]
Since \(a\in\Lambda=\operatorname{supp}\nu_A\) by
Lemma~\ref{model:factor}, we have \(\nu_A(U_a)>0\).
Using \(\Gamma_A(I)=A\) and \eqref{model:spectral-bound}, we obtain
\[
\begin{aligned}
  \|E(U_a)\Gamma_A(G_\Phi)\|_2
  &\ge\|E(U_a)A\|_2-\|E(U_a)\Gamma_A(G_\Phi-I)\|_2\\
  &\ge\nu_A(U_a)^{1/2}
      -\left(\int_{U_a}\|G_\Phi(w)-I\|_{\mathrm{op}}^2\,d\nu_A(w)\right)^{1/2}\\
  &\ge\tfrac14\nu_A(U_a)^{1/2}>0.
\end{aligned}
\]
Thus \(\Gamma_A(G_\Phi)\ne0\), so we may choose
\(v\in\mathcal E\) with \(\Gamma_A(G_\Phi)v\ne0\).
By \eqref{model:spectral-identities},
\[
  x:=\Gamma_A(G_\Phi)v
  =\Gamma_A(G_{\Phi v})
  =J(\Phi v)\ne0.
\]

Since \(\Delta Q=L\) almost everywhere on \(G\), applying
\eqref{cau:represent} to \((z-a)p(z)\) gives
\[
  \int_U p(z)\Delta(z)\Phi(z)\,dm_2(z)
  =L\int_G(z-a)p(z)h(z)\,dm_2(z)=0,
  \qquad p\in\mathbb C[z].
\]
For \(p\in\mathbb C[z]\), \eqref{model:intertwine} and the moment
identity give
\[
  TJ(p\Phi v)
  =J(zp\Phi v)
   -A\left(\int_U p(z)\Delta(z)\Phi(z)\,dm_2(z)\right)v
  =J(zp\Phi v).
\]
Since \(x=J(\Phi v)\), induction gives
\[
  T^n x=J(z^n\Phi v)\qquad(n\ge0).
\]
By linearity, \(p(T)x=J(p\Phi v)\) for every \(p\in\mathbb C[z]\).
By \eqref{model:averaging-bound},
\[
\begin{aligned}
  \|p(T)x\|
  =\|J(p\Phi v)\|
   \le c_U\|A\|_2\|p\Phi v\|_\infty
  \le c_U\|A\|_2\|\Phi\|_\infty\|v\|
       \sup_{z\in U}|p(z)|.
\end{aligned}
\]
Thus the assertion holds with
\(C_x=c_U\|A\|_2\|\Phi\|_\infty\|v\|<\infty\).
\end{proof}

Let $T=S+K$, where $S$ is subnormal on $H$ and $K\in S_1(H)$.
Let $N$ be the minimal normal extension of $S$ on $\widetilde H$.
Write $\iota:H\to\widetilde H$ for the inclusion and $P=\iota^*$
for the orthogonal projection onto $H$. Then $N\iota=\iota S$.

Since $K^*\in S_1(H)$, Lemma~\ref{model:factor}\textup{(1)} gives
Hilbert--Schmidt operators $A_0:\mathcal E\to H$ and
$B_0:H\to\mathcal E$ such that $K^*=A_0B_0$. \textcolor{black}{Set}
\[
  A=\iota A_0,\qquad B=B_0P.
\]
Set $W=N^*+AB=N^*+\iota K^*P$.
Let $E$ be the spectral measure of $N^*$, and \textcolor{black}{set} $\mu(L)=BE(L)A$
for every Borel set $L\subset\mathbb C$. By the proof of
Lemma~\ref{model:factor}\textup{(3)}, $\mu$ is countably additive in
$S_2(\mathcal E)$, and
\[
  |\mu|(\mathbb C)
  =\sup_{(L_j)}\sum_j\|BE(L_j)A\|_2
  \le\|B\|_2\|A\|_2,
\]
where the supremum is taken over all finite Borel partitions of
$\mathbb C$. Let $D$ be the principal-value Cauchy transform of $\mu$
defined in Lemma~\ref{model:regularized-transform}, and \textcolor{black}{set}
$\Delta(z)=I_{\mathcal E}-D(z)$.

\begin{lemma}\label{model:subnormal-extension}
With the notation above, the following hold.
\begin{enumerate}
\item
For every $p\in\mathbb C[z]$,
\begin{equation}\label{sub:extension}
  PW=T^*P,\qquad Pp(W)=p(T^*)P.
\end{equation}
\item
For each bounded open disk $U$ there is a bounded linear map
$J:L^\infty(U,m_2;\mathcal E)\to\widetilde H$ such that, for every $b$,
\begin{equation}\label{sub:transfer-bound}
  \|Jb\|\le c_U\|A\|_2\|b\|_\infty,\qquad
  BJb=\int_U D(z)b(z)\,dm_2(z),
\end{equation}
and
\begin{equation}\label{sub:transfer}
  WJb=J(zb)-A\int_U\Delta(z)b(z)\,dm_2(z).
\end{equation}
\end{enumerate}
\end{lemma}

\begin{proof}
\textup{(1)} Since $P=\iota^*$, we have
$PN^*=(N\iota)^*=(\iota S)^*=S^*P$. Thus $P\iota=I_H$ and
$A_0B_0=K^*$ give
\[
  PW=PN^*+(P\iota)A_0B_0P=(S^*+K^*)P=T^*P.
\]
By induction, $PW^j=(T^*)^jP$ for every $j\ge0$. Hence, for
$p(z)=\sum_{j=0}^m a_jz^j$,
\[
  Pp(W)=\sum_{j=0}^m a_jPW^j
       =\sum_{j=0}^m a_j(T^*)^jP
       =p(T^*)P,
\]
which proves \eqref{sub:extension}.

\smallskip
\noindent\textup{(2)} For $b\in L^\infty(U,m_2;\mathcal E)$, \textcolor{black}{set}
\[
  G_b(w)=\int_U\frac{b(z)}{z-w}\,dm_2(z)
  \qquad (w\in\mathbb C).
\]
By Lemma~\ref{model:averaging}\textup{(1)}, $G_b$ is norm-continuous.
We may therefore define $Jb=\Gamma_A(G_b|_{\sigma(N^*)})$, where
$\Gamma_A$ is formed with the spectral measure of $N^*$. With $c_U$
as in \eqref{cau:area}, the spectral estimate \eqref{model:spectral-bound}
gives
\[
\begin{aligned}
  \|Jb\|
  &\le \|A\|_2\sup_{w\in\sigma(N^*)}\|G_b(w)\| \\
  &\le \|A\|_2\|b\|_\infty
       \sup_{w\in\sigma(N^*)}\int_U\frac{dm_2(z)}{|z-w|}
   \le c_U\|A\|_2\|b\|_\infty.
\end{aligned}
\]

To compute $BJb$, note that \eqref{cau:area} also gives
\[
  \int_{\sigma(N^*)}\int_U
  \frac{\|b(z)\|}{|z-w|}\,dm_2(z)\,d|\mu|(w)
  \le c_U\|b\|_\infty|\mu|(\sigma(N^*))<\infty.
\]
Thus, as in the proof of Lemma~\ref{model:averaging}\textup{(2)},
\eqref{model:spectral-identities} and Fubini's theorem give
\[
\begin{aligned}
  BJb
  &=\int_{\sigma(N^*)}d\mu(w)\,G_b(w)\\
  &=\int_U\left(\int_{\sigma(N^*)}\frac{d\mu(w)}{z-w}\right)
      b(z)\,dm_2(z)\\
  &=\int_U D(z)b(z)\,dm_2(z).
\end{aligned}
\]
The last equality follows from Lemma~\ref{cau:area-lemma}, since the
inner integral equals $D(z)$ for $m_2$-almost every $z\in U$.

Now let $zb$ denote the field $z\mapsto z\,b(z)$. For $w\in\sigma(N^*)$,
\[
  G_{zb}(w)-wG_b(w)
  =\int_U\frac{(z-w)b(z)}{z-w}\,dm_2(z)
  =\int_U b(z)\,dm_2(z).
\]
Thus $G_{zb}-wG_b$ is constant on $\sigma(N^*)$. Applying
Lemma~\ref{model:spectral-integral}\textup{(2), (3)} gives
\[
  J(zb)-N^*Jb
  =\Gamma_A(G_{zb}-wG_b)
  =A\int_U b(z)\,dm_2(z).
\]
Finally, since $W=N^*+AB$ and $\Delta(z)=I_{\mathcal E}-D(z)$,
the formula for $BJb$ in \eqref{sub:transfer-bound} gives
\[
\begin{aligned}
  J(zb)-WJb
  &=A\int_U b(z)\,dm_2(z)-ABJb\\
  &=A\int_U\bigl(I_{\mathcal E}-D(z)\bigr)b(z)\,dm_2(z)\\
  &=A\int_U\Delta(z)b(z)\,dm_2(z),
\end{aligned}
\]
which proves \eqref{sub:transfer}.
\end{proof}

For the next lemma, let \(R\) be a bounded normal operator on
\(\widetilde H\), and let \(B\colon\widetilde H\to\mathcal E\) be
compact. Fix \(a\in\mathbb C\) and an orthonormal sequence \((s_n)\)
in \(\widetilde H\).

\begin{lemma}\label{model:spectral-localization}
Suppose that \(\|(R-aI)s_n\|\to0\). Let \(G\subset\mathbb C\) be a
nonempty Borel set with \(\operatorname{dist}(a,G)>0\). Then there
exist vectors \(y_n\in\widetilde H\) and, for all sufficiently
large \(n\), bounded continuous maps \(u_n\colon G\to\widetilde H\)
with the following properties.
\begin{enumerate}
\item
\(\displaystyle \lim_{n\to\infty}\|y_n-s_n\|=0\).

\item
For every \(z\in G\), \((zI-R)u_n(z)=y_n\), and
\[
  \lim_{n\to\infty}\sup_{z\in G}
  \left\|u_n(z)-\frac{y_n}{z-a}\right\|=0.
\]

\item
\(\displaystyle \lim_{n\to\infty}\sup_{z\in G}\|Bu_n(z)\|=0\).
\end{enumerate}
\end{lemma}

\begin{proof}
\textup{(1)} \textcolor{black}{Set} \(\varepsilon_n=(\|(R-aI)s_n\|+n^{-2})^{1/2}\). Since
\(\|(R-aI)s_n\|\to0\), we have
\[
  0\le\frac{\|(R-aI)s_n\|}{\varepsilon_n}
  \le\varepsilon_n\longrightarrow0.
\]
Let \(E\) be the spectral measure of \(R\), and \textcolor{black}{set}
\(y_n=E(B(a,\varepsilon_n))s_n\). By the spectral theorem,
\[
\begin{aligned}
  \|s_n-y_n\|^2
  &=\int_{\{|w-a|\ge\varepsilon_n\}}
      d\langle E(w)s_n,s_n\rangle\\
  &\le \frac{1}{\varepsilon_n^2}
       \int_{\mathbb C}|w-a|^2\,d\langle E(w)s_n,s_n\rangle\\
  &=\frac{\|(R-aI)s_n\|^2}{\varepsilon_n^2}
    \longrightarrow0.
\end{aligned}
\]

\smallskip
\noindent\textup{(2)} \textcolor{black}{Set} \(d=\operatorname{dist}(a,G)\). Since \(d>0\) and
\(\varepsilon_n\to0\), we have \(\varepsilon_n<d/2\) for all
sufficiently large \(n\). Since \(|z-a|\ge d\) on \(G\), the triangle inequality gives
\[
  |z-w|\ge |z-a|-|w-a|\ge d-\varepsilon_n>d/2,
  \qquad z\in G,\ w\in B(a,\varepsilon_n).
\]
Define
\[
  u_n(z)=\int_{B(a,\varepsilon_n)}\frac{dE(w)y_n}{z-w}
  \qquad (z\in G).
\]
By the spectral theorem, \(\|u_n(z)\|\le\|y_n\|/(d-\varepsilon_n)\).
For \(z,z'\in G\), the same estimate applied to the difference gives
\[
\begin{aligned}
  \|u_n(z)-u_n(z')\|
  &=\left\|\int_{B(a,\varepsilon_n)}
      \frac{z'-z}{(z-w)(z'-w)}\,dE(w)y_n\right\|\\
  &\le\frac{|z-z'|\,\|y_n\|}{(d-\varepsilon_n)^2}.
\end{aligned}
\]
Thus \(u_n\) is bounded and continuous on \(G\). Since
\(E(B(a,\varepsilon_n))y_n=y_n\), the spectral calculus gives
\[
  (zI-R)u_n(z)
  =\int_{B(a,\varepsilon_n)}\frac{z-w}{z-w}\,dE(w)y_n
  =y_n.
\]
Finally, subtracting \((z-a)^{-1}y_n\) and using
\(\|y_n\|\le\|s_n\|=1\), we obtain
\[
\begin{aligned}
  \sup_{z\in G}\left\|u_n(z)-\frac{y_n}{z-a}\right\|
  &=\sup_{z\in G}\left\|\int_{B(a,\varepsilon_n)}
      \frac{w-a}{(z-w)(z-a)}\,dE(w)y_n\right\|\\
  &\le\frac{\varepsilon_n\|y_n\|}{d(d-\varepsilon_n)}
   \longrightarrow0.
\end{aligned}
\]

\smallskip
\noindent\textup{(3)} Since \((s_n)\) is orthonormal,
\(s_n\rightharpoonup0\). By \textup{(1)}, \(\|y_n-s_n\|\to0\), so
\(y_n\rightharpoonup0\). The compactness of \(B\) gives
\(\|By_n\|\to0\). Hence, by \textup{(2)},
\[
  \sup_{z\in G}\|Bu_n(z)\|
  \le\|B\|\sup_{z\in G}
     \left\|u_n(z)-\frac{y_n}{z-a}\right\|
     +\frac{\|By_n\|}{d}
  \longrightarrow0.
\]
\end{proof}

\begin{lemma}\label{model:annular-weight}
Let \(\mu,D,\Delta\) be as in Lemma~\ref{model:subnormal-extension}.
Let \(U\subset\mathbb C\) be a bounded open disk and \(a\in U\).
Suppose that \(I_\mu(a)<\infty\) and that \(\Delta(a)\) is invertible.
Then there exist a Borel set \(G\subset U\) with
\(\operatorname{dist}(a,G)>0\) and a function
\(\varphi\in L^\infty(G,m_2)\) with the following properties.
\begin{enumerate}
\item
For \(m_2\)-almost every \(z\in G\), the inverse \(\Delta(z)^{-1}\)
exists and satisfies
\[
  \|\Delta(z)^{-1}\|_{\mathrm{op}}
  \le 2\|\Delta(a)^{-1}\|_{\mathrm{op}}.
\]
\item
The function \(\varphi\) satisfies
\[
  \int_G p(z)\varphi(z)\,dm_2(z)=0
  \quad (p\in\mathbb C[z]),\qquad
  \int_G \frac{\varphi(z)}{z-a}\,dm_2(z)=1.
\]
\end{enumerate}
\end{lemma}

\begin{proof}
Since \(\int_{\mathbb C}|a-w|^{-1}\,d|\mu|(w)=I_\mu(a)<\infty\),
dominated convergence gives
\[
  \lim_{t\to0^+}\mathcal C_t\mu(a)
  =\int_{\mathbb C}\frac{d\mu(w)}{a-w}
  =\mathcal C\mu(a)=D(a).
\]
By Lemma~\ref{cau:area-lemma}, \(D(z)=\mathcal C\mu(z)\) for
\(m_2\)-almost every \(z\in U\).
Choose \(0<\delta<(2\|\Delta(a)^{-1}\|_{\mathrm{op}})^{-1}\), and
let \(G\subset U\) be the Borel set given by
Lemma~\ref{cau:annular-density} for this \(\delta\). Then
\(\operatorname{dist}(a,G)>0\), and for \(m_2\)-almost every \(z\in G\),
\[
  \|D(z)-D(a)\|_2
  =\|\mathcal C\mu(z)-\mathcal C\mu(a)\|_2
  \le\delta.
\]
At each such \(z\), the identity \(\Delta=I_{\mathcal E}-D\) gives
\[
  \|I_{\mathcal E}-\Delta(a)^{-1}\Delta(z)\|_{\mathrm{op}}
  \le \|\Delta(a)^{-1}\|_{\mathrm{op}}\,\|D(z)-D(a)\|_2
  \le \delta\|\Delta(a)^{-1}\|_{\mathrm{op}}<\tfrac12.
\]
Hence \(\Delta(z)\) is invertible, and the Neumann series gives
\[
  \|\Delta(z)^{-1}\|_{\mathrm{op}}
  \le \frac{\|\Delta(a)^{-1}\|_{\mathrm{op}}}
           {1-\delta\|\Delta(a)^{-1}\|_{\mathrm{op}}}
  \le 2\|\Delta(a)^{-1}\|_{\mathrm{op}}.
\]
This proves \textup{(1)}.

For this same \(G\), Lemma~\ref{cau:annular-density} gives a function
\(h\in L^\infty(G,m_2)\) such that
\[
  \int_G p(z)h(z)\,dm_2(z)=p(a),
  \qquad p\in\mathbb C[z].
\]
\textcolor{black}{Set} \(\varphi(z)=(z-a)h(z)\) on \(G\). Since \(G\subset U\) and
\(U\) is bounded, \(\varphi\in L^\infty(G,m_2)\). Applying the
representing identity for \(h\) to the polynomial \((z-a)p(z)\), we get
\[
  \int_G p(z)\varphi(z)\,dm_2(z)
  =\int_G (z-a)p(z)h(z)\,dm_2(z)=0,
  \qquad p\in\mathbb C[z].
\]
Applying the same identity to the constant polynomial \(1\) gives
\[
  \int_G\frac{\varphi(z)}{z-a}\,dm_2(z)
  =\int_G h(z)\,dm_2(z)=1.
\]
This proves \textup{(2)}.
\end{proof}

\begin{lemma}\label{model:subnormal-regular}
Use the notation of Lemma~\ref{model:subnormal-extension}.
Let \(U\subset\mathbb C\) be a bounded open disk, and let \(a\in U\).
Suppose that the following conditions hold.
\begin{enumerate}
\item \(I_\mu(a)<\infty\).
\item \(\Delta(a)\) is invertible.
\item There is an orthonormal sequence \(\{s_n\}_{n=1}^{\infty}\) in \(H\) such that
\[
  \lim_{n\to\infty}\|(S-\bar a I)s_n\|=0.
\]
\end{enumerate}
Then there exist \(v\in H\setminus\{0\}\) and \(C_v<\infty\) such that
\[
  \|p(T^*)v\|\le C_v\sup_{z\in U}|p(z)|,
  \qquad p\in\mathbb C[z].
\]
\end{lemma}
\begin{proof}
By hypotheses \textup{(1)} and \textup{(2)},
Lemma~\ref{model:annular-weight} gives a Borel set \(G\subset U\)
with \(\operatorname{dist}(a,G)>0\).
We regard \(H\) as a closed subspace of \(\widetilde H\), on which
the minimal normal extension \(N\) acts. Since \(N-\bar a I\) is normal and \(N|_H=S\),
\[
  \|(N^*-aI)s_n\|
  =\|(N-\bar a I)^*s_n\|
  =\|(N-\bar a I)s_n\|
  =\|(S-\bar a I)s_n\|
  \longrightarrow0.
\]
Since \(B\) is Hilbert--Schmidt and hence compact,
Lemma~\ref{model:spectral-localization} applies to \(N^*\) and \((s_n)\). By parts \textup{(1)} and
\textup{(2)} of Lemma~\ref{model:spectral-localization}, there are
vectors \(y_n\in\widetilde H\) with \(\|y_n-s_n\|\to0\) and, for
each sufficiently large \(n\), a bounded continuous map
\(u_n\colon G\to\widetilde H\) such that
\[
  (zI-N^*)u_n(z)=y_n, \quad z\in G.
\]
Let \(\varphi\in L^\infty(G,m_2)\) be the function furnished by
Lemma~\ref{model:annular-weight} on \(G\). Since \(W=N^*+AB\) and
\((zI-N^*)u_n(z)=y_n\), we have, almost everywhere on \(G\),
\[
  (zI-W)\bigl(\varphi(z)u_n(z)\bigr)
  =\varphi(z)y_n-A\bigl(\varphi(z)Bu_n(z)\bigr).
\] For each such \(n\),
\textcolor{black}{let}
\[
  b_n(z)=\varphi(z)\Delta(z)^{-1}Bu_n(z)
  \qquad\text{for \(m_2\)-almost every \(z\in G\),}
\]
and extend \(b_n\) by zero on \(U\setminus G\). Since \(u_n\) is
bounded and \(\Delta^{-1}\) is essentially bounded on \(G\) by
Lemma~\ref{model:annular-weight}\textup{(1)}, we have
\(b_n\in L^\infty(U,m_2;\mathcal E)\).

Let \(J\) be the map in Lemma~\ref{model:subnormal-extension} for
the disk \(U\). For \(p\in\mathbb C[z]\), \textcolor{black}{set}
\[
  x_{n,p}=\int_G p(z)\varphi(z)u_n(z)\,dm_2(z)+J(pb_n).
\]
Applying \eqref{sub:transfer} to \(pb_n\) and using
\((zI-N^*)u_n(z)=y_n\), we obtain
\[
\begin{aligned}
  x_{n,zp}-Wx_{n,p}
  &=\int_G p(z)\varphi(z)(zI-W)u_n(z)\,dm_2(z)
    +A\int_G p\Delta b_n\,dm_2\\
  &=y_n\int_G p\varphi\,dm_2
    -A\int_G p\varphi Bu_n\,dm_2
    +A\int_G p\Delta b_n\,dm_2\\
  &=0.
\end{aligned}
\]
Indeed, \(\int_G p\varphi\,dm_2=0\) by
Lemma~\ref{model:annular-weight}\textup{(2)}, and the last two terms
cancel because \(\Delta b_n=\varphi Bu_n\) almost everywhere on \(G\).

It follows by induction that \(x_{n,z^k}=W^k x_{n,1}\) for every
\(k\ge0\). Since \(x_{n,p}\) depends linearly on \(p\),
\[
  x_{n,p}=p(W)x_{n,1},\qquad p\in\mathbb C[z].
\]
By the definition of \(x_{n,p}\),
\[
\begin{aligned}
  \|x_{n,p}\|
  &\le \int_G |p(z)|\,|\varphi(z)|\,\|u_n(z)\|\,dm_2(z)
       +\|J\|\,\|pb_n\|_\infty\\
  &\le C_n\sup_{z\in U}|p(z)|,
\end{aligned}
\]
where
\[
  C_n=\|\varphi\|_1\sup_{z\in G}\|u_n(z)\|
      +\|J\|\,\|b_n\|_\infty<\infty.
\]

By \eqref{sub:extension} and the identity \(x_{n,p}=p(W)x_{n,1}\),
\[
  p(T^*)Px_{n,1}=Pp(W)x_{n,1}=Px_{n,p},
  \qquad p\in\mathbb C[z].
\]
Since \(P\) is a contraction, the preceding estimate gives
\[
  \|p(T^*)Px_{n,1}\|\le\|x_{n,p}\|
  \le C_n\sup_{z\in U}|p(z)|.
\]
It remains to show that \(Px_{n,1}\ne0\) for some \(n\).

By Lemma~\ref{model:annular-weight}\textup{(2)},
\(\int_G\varphi(z)/(z-a)\,dm_2(z)=1\). Hence
\[
  x_{n,1}-s_n
  =\int_G\varphi(z)
     \left(u_n(z)-\frac{y_n}{z-a}\right)\,dm_2(z)
    +Jb_n+(y_n-s_n).
\]
To estimate \(Jb_n\), note that
Lemma~\ref{model:annular-weight}\textup{(1)} gives
\[
  \|b_n\|_\infty
  \le 2\|\varphi\|_\infty\|\Delta(a)^{-1}\|_{\mathrm{op}}
      \sup_{z\in G}\|Bu_n(z)\|.
\]
By Lemma~\ref{model:spectral-localization}\textup{(3)},
\(\lim\limits_{n\to\infty}\sup_{z\in G}\|Bu_n(z)\|=0\). It follows that
\(\lim\limits_{n\to\infty}\|b_n\|_\infty=0\). Taking norms in the identity
for \(x_{n,1}-s_n\) and applying \eqref{sub:transfer-bound}, we obtain
\[
  \|x_{n,1}-s_n\|
  \le \|\varphi\|_1
      \sup_{z\in G}\left\|u_n(z)-\frac{y_n}{z-a}\right\|
      +c_U\|A\|_2\|b_n\|_\infty+\|y_n-s_n\|
  \longrightarrow0,
\]
where the first and third terms tend to zero by
Lemma~\ref{model:spectral-localization}\textup{(2)} and
\textup{(1)}, respectively.
Since \(Ps_n=s_n\) and \(\|s_n\|=1\),
\[
  \bigl|\|Px_{n,1}\|-1\bigr|
  \le \|P(x_{n,1}-s_n)\|
  \le \|x_{n,1}-s_n\|\longrightarrow0.
\]
Fix \(n\) for which \(Px_{n,1}\ne0\), and \textcolor{black}{set} \(v=Px_{n,1}\).
Then the assertion holds with \(C_v=C_n\).
\end{proof}

\section{Localization on sets of positive analytic capacity}\label{sec:localization}

We first use the notation of Section~\ref{sec:model} for the normal
operator \(N\) and its trace-class perturbation \(T\).
Fix a bounded open disk \(U\) and a compact
set \(L_0\subset U\cap\Lambda\) with \(\gamma(L_0)>0\).
By Lemma~\ref{cau:capacity}\textup{(2)}, there is a finite positive
measure \(\nu\) supported on \(L_0\) satisfying \eqref{cau:measure}.
In particular, \(\nu(L_0)\ge c\gamma(L_0)>0\).

We shall use the following Hilbert-space-valued version of the
Davie--{\O}ksendal cutoff lemma. For the classical scalar cutoff theorem with \(f\equiv1\),
see \cite[Theorem~2.1]{TolsaICM}.

\begin{lemma}\label{gen:cutoff}
Let \(F\subset U\) be compact with \(\nu(F)>0\), and let
\(f\colon F\to\mathcal E\) be bounded and strongly measurable.
There exists a measurable function \(h\colon F\to[0,1]\) such that
\[
  \int_F h\,d\nu\ge\tfrac12\nu(F),\qquad
  \sup_{\substack{w\in\mathbb C\\ r>0}}
  \left\|\int_F k_r(z-w)h(z)f(z)\,d\nu(z)\right\|_{\mathcal E}
  <\infty.
\]
\end{lemma}
\begin{proof}
If \(\|f\|_\infty=0\), take \(h=1\). Assume that \(\|f\|_\infty>0\).
Fix \((w_j,r_j)\in\mathbb C\times(0,\infty)\), \(1\le j\le n\),
and let \(\beta_j\) be normalized area measure on \(B(w_j,r_j)\).
Then
\[
  \mathcal C\beta_j(z)
  =\frac{1}{\pi r_j^2}\int_{B(w_j,r_j)}
       \frac{dm_2(\zeta)}{z-\zeta}
  =k_{r_j}(z-w_j).
\]
For \(\ell=(\ell_1,\ldots,\ell_n)\in(\mathcal E^*)^n\) with
\(\sum\limits_{j=1}^{n}\|\ell_j\|\le1\), \textcolor{black}{set}
\(\eta=\sum\limits_{j=1}^{n}\ell_j\beta_j\).
Since \(\eta\) has a bounded, compactly supported density with respect
to \(m_2\) and \(|\cdot|^{-1}\in L^1_{\mathrm{loc}}(\mathbb C,m_2)\),
the integral defining \(\mathcal C\eta(z)\) converges absolutely for
every \(z\in\mathbb C\). For \(z\in F\), \textcolor{black}{set}
\[
\begin{aligned}
  s_\ell(z)
  &:=\sum\limits_{j=1}^{n} k_{r_j}(z-w_j)\ell_j(f(z))\\
  &=\left(\sum\limits_{j=1}^{n}\mathcal C\beta_j(z)\,\ell_j\right)(f(z))
   =\bigl(\mathcal C\eta(z)\bigr)(f(z)).
\end{aligned}
\]
Hence, for \(\nu\)-almost every \(z\in F\),
\[
  |s_\ell(z)|
  =\bigl|(\mathcal C\eta(z))(f(z))\bigr|
  \le \|\mathcal C\eta(z)\|_{\mathcal E^*}\|f(z)\|
  \le \|f\|_\infty\,\mathcal C_*\eta(z).
\]
Since \(\beta_j(\mathbb C)=1\), we have
\(|\eta|(\mathbb C)\le\sum\limits_{j=1}^{n}\|\ell_j\|\le1\).
For \(t>0\), \textcolor{black}{let}
\(E_t=\{z\in\mathbb C:\mathcal C_*\eta(z)>t/\|f\|_\infty\}\). Then
\begin{equation}\label{gen:cutoff-weak}
\begin{aligned}
  \nu\{z\in F:|s_\ell(z)|>t\}
  &\le \nu(E_t)
    \overset{\text{\scriptsize\eqref{cau:borel-domination}}}{\le} C_1\gamma(E_t)\\
  &\mathrel{\overset{\text{\scriptsize Lemma~\ref{cau:max}}}{\le}}
    \frac{C\|f\|_\infty}{t}\,|\eta|(\mathbb C)
   \le \frac{C\|f\|_\infty}{t},
\end{aligned}
\end{equation}
where \(C\) is independent of \(\ell\) and of the finite family.
Take \(t=2C\|f\|_\infty/\nu(F)\) and \textcolor{black}{let}
\(h_\ell=\mathbf1_{\{|s_\ell|\le t\}}\) on \(F\).
By \eqref{gen:cutoff-weak},
\[
  \int_F h_\ell\,d\nu
  =\nu(F)-\nu\{z\in F:|s_\ell(z)|>t\}
  \ge\tfrac12\nu(F).
\]
Moreover,
\[
  \operatorname{Re}\int_F h_\ell s_\ell\,d\nu
  \le \int_F h_\ell|s_\ell|\,d\nu
  \le t\,\nu(F)=2C\|f\|_\infty.
\]

We now apply the Davie--{\O}ksendal duality argument
\cite[Section~4.3]{TolsaBook}. \textcolor{black}{Set}
\[
  \mathcal K=
  \left\{h\in L^\infty(\nu|_F):0\le h\le1,\quad
         \int_F h\,d\nu\ge\tfrac12\nu(F)\right\}.
\]
Then \(h_\ell\in\mathcal K\) for every \(\ell\) as above.
The set \(\mathcal K\) is weak-star compact and convex, and the pairing
\[
  (h,\ell)\longmapsto\operatorname{Re}\int_F h s_\ell\,d\nu
\]
is real bilinear, weak-star continuous in \(h\), and norm-continuous
in \(\ell\). By Sion's minimax theorem~\textcolor{black}{\cite[p.~174]{Sion}}
and the preceding bound, we obtain
\[
\begin{aligned}
  \min_{h\in\mathcal K}\sup_\ell
    \operatorname{Re}\int_F h s_\ell\,d\nu
  &= \sup_\ell\min_{h\in\mathcal K}
    \operatorname{Re}\int_F h s_\ell\,d\nu\\
  &\le \sup_\ell \operatorname{Re}\int_F h_\ell s_\ell\,d\nu
  \le 2C\|f\|_\infty,
\end{aligned}
\]
where each supremum is taken over all \(\ell=(\ell_1,\ldots,\ell_n)\)
satisfying \(\sum\limits_{j=1}^n\|\ell_j\|\le 1\).
For a minimizing \(h\in\mathcal K\), duality gives
\begin{equation}\label{gen:cutoff-finite}
  \max_{1\le j\le n}
  \left\|\int_F k_{r_j}(z-w_j)h(z)f(z)\,d\nu(z)\right\|_{\mathcal E}
  =\sup_\ell\operatorname{Re}\int_F h s_\ell\,d\nu
  \le 2C\|f\|_\infty.
\end{equation}

For \(w\in\mathbb C\) and \(r>0\), \textcolor{black}{set}
\[
  \mathcal K_{w,r}
  =\left\{h\in\mathcal K:
    \left\|\int_F k_r(z-w)h(z)f(z)\,d\nu(z)\right\|_{\mathcal E}
    \le 2C\|f\|_\infty\right\}.
\]
By \eqref{gen:cutoff-finite},
\[
  \bigcap_{j=1}^n\mathcal K_{w_j,r_j}\ne\varnothing
\]
for every finite family \((w_j,r_j)_{j=1}^n\). By duality,
\[
  \left\|\int_F k_r(z-w)h(z)f(z)\,d\nu(z)\right\|_{\mathcal E}
  =
  \sup_{\substack{\ell\in\mathcal E^*\\ \|\ell\|\le1}}
  \left|\int_F h(z)k_r(z-w)\ell(f(z))\,d\nu(z)\right|.
\]
Since \(k_r(\cdot-w)\ell(f)\in L^1(\nu|_F)\), each
\(\mathcal K_{w,r}\) is weak-star closed in \(\mathcal K\).
By weak-star compactness of \(\mathcal K\),
\[
  \bigcap_{\substack{w\in\mathbb C\\ r>0}}\mathcal K_{w,r}
  \ne\varnothing.
\]
Choose an element of this intersection and let \(h\colon F\to[0,1]\)
be a measurable representative. Since \(h\in\mathcal K\), we have
\(0\le h\le1\) and \(\int_F h\,d\nu\ge\tfrac12\nu(F)\).
Since \(h\in\mathcal K_{w,r}\) for every \(w\in\mathbb C\) and \(r>0\),
\[
  \sup_{\substack{w\in\mathbb C\\ r>0}}
  \left\|\int_F k_r(z-w)h(z)f(z)\,d\nu(z)\right\|_{\mathcal E}
  \le 2C\|f\|_\infty.
\]
\end{proof}
Let \(\mu,D,\Delta\) be as in Section~\ref{sec:model}.

\begin{lemma}\label{gen:polynomial-moments}
Let \(F\subset U\cap\Lambda\) be compact with \(\nu(F)>0\), let
\(g\in L^\infty(\nu|_F;\mathcal E)\), and let \(r_n\rightarrow0\).
Suppose that
\[
\sup_{z\in F}\|\widetilde{\mathcal C}_{r_n}\mu(z)-D(z)\|_2
\longrightarrow0,\qquad
\sup_{n\ge1,\,w\in\mathbb C}
\left\|\int_F k_{r_n}(z-w)g(z)\,d\nu(z)\right\|<\infty.
\]
If \(u\in E(F)H\) satisfies
\begin{equation}\label{gen:correction}
Bp(N)u=\int_F p(z)\Delta(z)g(z)\,d\nu(z),
\qquad p\in\mathbb C[z],
\end{equation}
then there exist \(x\in H\) and \(C_x<\infty\) such that, for every
\(p\in\mathbb C[z]\),
\[
\|p(T)x\|\le C_x\sup_{z\in U}|p(z)|,\qquad
Bp(T)x=\int_F p(z)g(z)\,d\nu(z).
\]
\end{lemma}

\begin{proof}
For \(p\in\mathbb C[z]\) and \(n\ge1\), \textcolor{black}{set}
\[
G_{n,p}(w)=\int_F k_{r_n}(z-w)p(z)g(z)\,d\nu(z).
\]
Since \(k_{r_n}\) is bounded and continuous, dominated convergence
shows that \(G_{n,p}\) is norm-continuous.
Choose a closed disk \(L\subset U\) with \(F\subset\operatorname{int}L\).
If \(z\in F\) and \(w\in L\), the segment from \(w\) to \(z\) lies
in \(L\), so that
\[
|p(z)-p(w)|\le |z-w|\sup_L|p'|.
\]
Write \(\|p\|_U=\sup\limits_{z\in U}|p(z)|\). Since
\(|s k_{r_n}(s)|\le1\), Cauchy's estimate gives
\[
\begin{aligned}
\|G_{n,p}(w)-p(w)G_{n,1}(w)\|
&\le \int_F |k_{r_n}(z-w)|\,|p(z)-p(w)|\,\|g(z)\|\,d\nu(z)\\
&\le \|g\|_{L^1(\nu|_F)}\sup_L|p'|\\
&\le \frac{\|g\|_{L^1(\nu|_F)}}{\operatorname{dist}(L,\partial U)}
\|p\|_U.
\end{aligned}
\]
Thus, for \(w\in L\),
\[
\|G_{n,p}(w)\|
\le \left(\|G_{n,1}\|_\infty+
\frac{\|g\|_{L^1(\nu|_F)}}{\operatorname{dist}(L,\partial U)}
\right)\|p\|_U.
\]
If \(w\notin L\), then
\(|z-w|\ge\operatorname{dist}(F,\mathbb C\setminus L)>0\)
for every \(z\in F\). Since \(|k_{r_n}(z-w)|\le |z-w|^{-1}\),
\[
\|G_{n,p}(w)\|
\le \frac{\|g\|_{L^1(\nu|_F)}}
{\operatorname{dist}(F,\mathbb C\setminus L)}
\|p\|_U.
\]
The two \textcolor{black}{bounds} cover all \(w\in\mathbb C\). Since
\(\sup\limits_n\|G_{n,1}\|_\infty<\infty\) by hypothesis, it follows that
\[
\sup_n\|G_{n,p}\|_\infty\le C_0\|p\|_U,
\]
where \(C_0\) is independent of \(p\).
\textcolor{black}{Let} \(X_n(p)=p(N)u+\Gamma_A(G_{n,p}|_\Lambda)\).
By \eqref{model:spectral-bound} and the identity
\(\nu_A(\Lambda)=\|A\|_2^2\) in Lemma~\ref{model:factor}\textup{(2)},
\[
\|\Gamma_A(G_{n,p}|_\Lambda)\|
\le \nu_A(\Lambda)^{1/2}\sup_{w\in\Lambda}\|G_{n,p}(w)\|
\le C_0\|A\|_2\|p\|_U.
\]
Also, \(E(F)u=u\) and \(F\subset U\) give
\(\|p(N)u\|\le\|u\|\|p\|_U\). Hence
\begin{equation}\label{gen:bounded}
\sup_n\|X_n(p)\|
\le \bigl(\|u\|+C_0\|A\|_2\bigr)\|p\|_U.
\end{equation}
Taking \(p=1\) in \eqref{gen:bounded}, we obtain
\[
\sup_n\|X_n(1)\|\le \|u\|+C_0\|A\|_2.
\]
Choose a subsequence such that \(X_{n_k}(1)\rightharpoonup x\) in \(H\).
To show that \(X_{n_k}(p)\rightharpoonup p(T)x\) for every polynomial
\(p\), we first estimate \(X_n(zp)-TX_n(p)\).

For each fixed \(n\), \(\lvert k_{r_n}(z-w)\rvert\le r_n^{-1}\). Hence
\[
  \int_\Lambda\int_F
    |k_{r_n}(z-w)|\,|p(z)|\,\|g(z)\|\,d\nu(z)\,d|\mu|(w)
  \le \frac{|\mu|(\Lambda)}{r_n}\,
       \|p\|_U\|g\|_{L^1(\nu|_F)}<\infty.
\]
By \eqref{model:spectral-identities} and Fubini's theorem,
\[
\begin{aligned}
  B\Gamma_A(G_{n,p}|_\Lambda)
  &=\int_\Lambda d\mu(w)\,G_{n,p}(w)\\
  &=\int_F p(z)
       \left(\int_\Lambda k_{r_n}(z-w)\,d\mu(w)\right)
       g(z)\,d\nu(z)\\
  &=\int_F p(z)\bigl(\widetilde{\mathcal C}_{r_n}\mu(z)\bigr)
       g(z)\,d\nu(z).
\end{aligned}
\]
\textcolor{black}{Let} \(c_p=\int_F p(z)g(z)\,d\nu(z)\). Since
\(BX_n(p)=Bp(N)u+B\Gamma_A(G_{n,p}|_\Lambda)\),
\eqref{gen:correction} and the identity \(\Delta=I_{\mathcal E}-D\)
give
\[
\begin{aligned}
  BX_n(p)-c_p
  &=\int_F p(z)\bigl(\Delta(z)+\widetilde{\mathcal C}_{r_n}\mu(z)
       -I_{\mathcal E}\bigr)g(z)\,d\nu(z)\\
  &=\int_F p(z)\bigl(\widetilde{\mathcal C}_{r_n}\mu(z)-D(z)\bigr)
       g(z)\,d\nu(z).
\end{aligned}
\]
By the assumed uniform convergence on \(F\),
\begin{equation}\label{gen:moment-limit}
  \|BX_n(p)-c_p\|
  \le \|p\|_U\|g\|_{L^1(\nu|_F)}
       \sup_{z\in F}\|\widetilde{\mathcal C}_{r_n}\mu(z)-D(z)\|_2
  \longrightarrow 0.
\end{equation}
To compare \(X_n(zp)\) with \(TX_n(p)\), \textcolor{black}{set}
\[
\begin{aligned}
  R_{n,p}(w)
  &=G_{n,zp}(w)-wG_{n,p}(w)-c_p\\
  &=\int_F \bigl((z-w)k_{r_n}(z-w)-1\bigr)p(z)g(z)\,d\nu(z).
\end{aligned}
\]
By \eqref{model:regularization},
\[
  \bigl|(z-w)k_{r_n}(z-w)-1\bigr|
  \le \mathbf{1}_{B(w,r_n)}(z).
\]
The growth bound in \eqref{cau:measure} therefore gives
\[
  \|R_{n,p}(w)\|
  \le \|p\|_U\|g\|_\infty\,\nu(B(w,r_n))
  \le r_n\|p\|_U\|g\|_\infty.
\]
In particular, \(\|R_{n,p}\|_\infty\to0\). Since \(T=N+AB\),
\eqref{model:spectral-identities} and the identity
\(\Gamma_A(c_p)=Ac_p\) in
Lemma~\ref{model:spectral-integral}\textup{(3)} give
\[
\begin{aligned}
  X_n(zp)-TX_n(p)
  &=\Gamma_A(R_{n,p}|_\Lambda)+Ac_p-ABX_n(p)\\
  &=\Gamma_A(R_{n,p}|_\Lambda)+A(c_p-BX_n(p)).
\end{aligned}
\]
It follows from \eqref{model:spectral-bound} and
\eqref{gen:moment-limit} that
\begin{equation}\label{gen:recurrence-limit}
  \|X_n(zp)-TX_n(p)\|
  \le \|A\|_2\|R_{n,p}\|_\infty
       +\|A\|\|c_p-BX_n(p)\|
  \longrightarrow 0.
\end{equation}

Recall that \(X_{n_k}(1)\rightharpoonup x\). Taking \(p=z^j\) in
\eqref{gen:recurrence-limit} and passing to the chosen subsequence, we obtain
\[
  \lim_{k\to\infty}
  \bigl\|X_{n_k}(z^{j+1})-TX_{n_k}(z^j)\bigr\|=0,
  \qquad j\ge0.
\]
Suppose that \(X_{n_k}(z^j)\rightharpoonup T^jx\). Since \(T\) is
bounded, \(TX_{n_k}(z^j)\rightharpoonup T^{j+1}x\), and hence
\[
  X_{n_k}(z^{j+1})
  =TX_{n_k}(z^j)
   +\bigl(X_{n_k}(z^{j+1})-TX_{n_k}(z^j)\bigr)
  \rightharpoonup T^{j+1}x.
\]
Induction gives \(X_{n_k}(z^j)\rightharpoonup T^jx\) for every \(j\ge0\).
By linearity,
\begin{equation}\label{gen:polynomial-weak-limit}
  X_{n_k}(p)\rightharpoonup p(T)x,
  \qquad p\in\mathbb C[z].
\end{equation}
Applying weak lower semicontinuity of the norm to
\eqref{gen:polynomial-weak-limit} and then using \eqref{gen:bounded}, we get
\[
  \|p(T)x\|
  \le \liminf_{k\to\infty}\|X_{n_k}(p)\|
  \le \bigl(\|u\|+C_0\|A\|_2\bigr)\|p\|_U.
\]

Finally, for every \(y\in\mathcal E\), \eqref{gen:moment-limit} and
\eqref{gen:polynomial-weak-limit} give
\[
\begin{aligned}
  \langle c_p,y\rangle
  &=\lim_{k\to\infty}\langle BX_{n_k}(p),y\rangle\\
  &=\lim_{k\to\infty}\langle X_{n_k}(p),B^*y\rangle\\
  &=\langle p(T)x,B^*y\rangle
   =\langle Bp(T)x,y\rangle.
\end{aligned}
\]
Thus
\[
  Bp(T)x=c_p=\int_F p(z)g(z)\,d\nu(z),
  \qquad p\in\mathbb C[z].
\]
\end{proof}

Write the Lebesgue decomposition of \(\mu\) with respect to \(\nu\) as
\[
d\mu=\mathcal M\,d\nu+d\mu_s,\qquad
\mathcal M\in L^1(\nu;S_2(\mathcal E)),\qquad
|\mu_s|\perp\nu,
\]
and choose a Borel representative of \(\mathcal M\).
By Lemma~\ref{cau:pv} and the definition of \(D\) in
Lemma~\ref{model:regularized-transform}, there is a Borel set \(Z\subset\mathbb C\)
with \(\gamma(Z)=0\) such that
\[
\lim_{t\downarrow0}\|\mathcal C_t\mu(z)-D(z)\|_2=0,
\qquad z\notin Z.
\]
By \eqref{cau:borel-domination}, \(\nu(Z)\le C\gamma(Z)=0\). \textcolor{black}{Set}
\[
\mathcal S=\{z\in\mathbb C\setminus Z:\Delta(z)\text{ is singular}\}
\]
and
\[
\mathcal V=\{z\in\mathbb C\setminus Z:
\Delta(z)\text{ is invertible and }\mathcal M(z)\ne0\}.
\]
Both sets are Borel, since \(\Delta\) and \(\mathcal M\) are norm-Borel
and the invertible operators form an open set.

\begin{lemma}\label{gen:compact-selection}
Suppose that \(\nu(\mathcal S\cup\mathcal V)>0\). Then there is
a compact set \(F\subset L_0\cap(\mathcal S\cup\mathcal V)\) such that
\[
\nu(F)>0,\qquad |\mu_s|(F)=0,\qquad
\lim_{n\to\infty}\sup_{z\in F}
\|\widetilde{\mathcal C}_{1/n}\mu(z)-D(z)\|_2=0.
\]
Moreover, there exist vectors \(e,v\in\mathcal E\) with \(\|v\|=1\),
a constant \(d>0\), and a bounded Borel map \(f\colon F\to\mathcal E\)
such that
\[
\Delta(z)f(z)=\mathcal M(z)e,\qquad
\operatorname{Re}\langle f(z),v\rangle\ge d,
\qquad z\in F.
\]
\end{lemma}

\begin{proof}
Since \(\nu(\mathcal S\cup\mathcal V)>0\), choose
\(W\in\{\mathcal S,\mathcal V\}\) with \(\nu(W)>0\).
We first show that there exist a fixed vector \(e\in\mathcal E\)
and a Borel map \(f\colon W\to\mathcal E\) such that
\[
\Delta(z)f(z)=\mathcal M(z)e,
~ z\in W,\qquad
\nu\{z\in W:f(z)\ne0\}>0.
\]
Fix an orthonormal basis \(\{b_j\}_{j=1}^{\infty}\) of \(\mathcal E\).
For each basis vector \(b_j\), we construct a vector
\(e_j\in\mathcal E\) and a Borel map \(f_j\colon W\to\mathcal E\).

\smallskip
\noindent\textbf{Case 1.} \(W=\mathcal S\).
Fix \(z\in W\). By Lemma~\ref{model:regularized-transform}\textup{(1)},
\(\Delta(z)\) is Fredholm of index zero. Hence
\[
  \dim\ker\Delta(z)
  =\operatorname{codim}\operatorname{ran}\Delta(z).
\]
If \(\ker\Delta(z)=\{0\}\), then \(\operatorname{ran}\Delta(z)=\mathcal E\),
and the bounded inverse theorem gives that \(\Delta(z)\) is invertible.
Since \(z\in\mathcal S\), this contradicts the singularity of \(\Delta(z)\).
Thus \(\ker\Delta(z)\ne\{0\}\).
Let \(P(z)\) be the orthogonal projection of \(\mathcal E\) onto
\(\ker\Delta(z)\). Since \(\ker(\Delta(z)^*\Delta(z))=\ker\Delta(z)\),
the spectral theorem gives
\[
  P(z)=\mathbf{1}_{\{0\}}\bigl(\Delta(z)^*\Delta(z)\bigr).
\]
For every \(t\ge0\), we have \(0<(1+kt)^{-1}\le1\) and
\((1+kt)^{-1}\to\mathbf{1}_{\{0\}}(t)\) as \(k\to\infty\).
Thus, by bounded convergence in the spectral calculus,
\[
  P(z)b
  =\lim_{k\to\infty}
   \bigl(I_{\mathcal E}+k\Delta(z)^*\Delta(z)\bigr)^{-1}b,
  \qquad b\in\mathcal E.
\]
For fixed \(k\) and \(b\), the map
\(A\mapsto(I_{\mathcal E}+kA^*A)^{-1}b\) is norm-continuous.
Since \(\Delta\) is norm-Borel, each term on the right is Borel in \(z\).
Thus \(z\mapsto P(z)b\), being the pointwise norm limit of these maps,
is Borel on \(W\).
For each \(j\ge1\), \textcolor{black}{let} \(f_j(z)=P(z)b_j\) on \(W\).
Then \(f_j\) is Borel and \(f_j(z)\in\operatorname{ran}P(z)=\ker\Delta(z)\).
Taking \(e_j=0\), we obtain
\[
  \Delta(z)f_j(z)=\Delta(z)P(z)b_j=0=\mathcal M(z)e_j,
  \qquad z\in W.
\]

\smallskip
\noindent\textbf{Case 2.} \(W=\mathcal V\).
For \(z\in W\), \(\Delta(z)\) is invertible. For each \(j\ge1\), set
\(e_j=b_j\) and \textcolor{black}{set}
\[
  f_j(z)=\Delta(z)^{-1}\mathcal M(z)b_j,
  \qquad z\in W.
\]
Since \(\Delta\) and \(\mathcal M\) are norm-Borel and inversion is
continuous on the invertible operators, \(f_j\) is Borel by continuity
of operator multiplication and evaluation at \(b_j\). Moreover,
\[
  \Delta(z)f_j(z)
  =\Delta(z)\Delta(z)^{-1}\mathcal M(z)b_j
  =\mathcal M(z)e_j,
  \qquad z\in W.
\]

\smallskip
Fix \(z\in W\). In Case~1, \(\operatorname{ran}P(z)=\ker\Delta(z)\ne\{0\}\).
Since \(\{b_j\}_{j\ge1}\) is an orthonormal basis of \(\mathcal E\),
there is some \(j\) for which \(f_j(z)=P(z)b_j\ne0\).
In Case~2, \(\mathcal M(z)\ne0\), so \(\mathcal M(z)b_j\ne0\) for
some \(j\). For this \(j\), the invertibility of \(\Delta(z)\) gives
\[
  f_j(z)=\Delta(z)^{-1}\mathcal M(z)b_j\ne0.
\]
Thus every \(z\in W\) belongs to at least one of the sets
\(\{z\in W:f_j(z)\ne0\}\), and hence
\[
  W=\bigcup_{j\ge1}\{z\in W:f_j(z)\ne0\}.
\]
Since \(\nu\) is supported on \(L_0\), countable subadditivity gives
\[
  0<\nu(W)=\nu(W\cap L_0)
  \le \sum_{j\ge1}\nu\{z\in W\cap L_0:f_j(z)\ne0\}.
\]
Hence there is a fixed index \(j_0\ge1\) such that
\[
  \nu\{z\in W\cap L_0:f_{j_0}(z)\ne0\}>0.
\]
\textcolor{black}{Let} \(f=f_{j_0}\) and \(e=e_{j_0}\). Then \(f\) is Borel and nonzero
on a subset of \(W\cap L_0\) of positive \(\nu\)-measure. Moreover,
\[
  \Delta(z)f(z)
  =\Delta(z)f_{j_0}(z)
  =\mathcal M(z)e_{j_0}
  =\mathcal M(z)e,
  \qquad z\in W.
\]

To obtain a uniform lower bound in a fixed direction, let
\(\{v_i\}_{i\ge1}\) be dense in the unit sphere of \(\mathcal E\).
For \(i,k\ge1\), \textcolor{black}{set}
\[
F_{i,k}
=\bigl\{z\in W\cap L_0:\|f(z)\|\le k,\
\operatorname{Re}\langle f(z),v_i\rangle\ge k^{-1}\bigr\}.
\]
Since \(\sup\limits_i\operatorname{Re}\langle y,v_i\rangle=\|y\|\)
for every \(y\in\mathcal E\),
\[
\{z\in W\cap L_0:f(z)\ne0\}
=\bigcup_{i,k\ge1}F_{i,k}.
\]
By countable subadditivity,
\[
0<\nu\{z\in W\cap L_0:f(z)\ne0\}
\le \sum_{i,k\ge1}\nu(F_{i,k}).
\]
We may therefore fix \(i,k\ge1\) such that \(\nu(F_{i,k})>0\).
Since \(|\mu_s|\perp\nu\), there is a Borel set \(Z_s\) such that
\[
  \nu(Z_s)=0,\qquad |\mu_s|(\mathbb C\setminus Z_s)=0.
\]
\textcolor{black}{Set} \(F_0=F_{i,k}\setminus Z_s\). Then
\[
  \nu(F_0)=\nu(F_{i,k})>0,\qquad |\mu_s|(F_0)=0.
\]
Since \(F_0\subset W\subset\mathbb C\setminus Z\),
\eqref{model:regularized-limit} gives
\[
  \lim_{n\to\infty}
  \|\widetilde{\mathcal C}_{1/n}\mu(z)-D(z)\|_2=0,
  \qquad z\in F_0.
\]
As \(0<\nu(F_0)<\infty\), Egorov's theorem gives a measurable set
\(F_1\subset F_0\) such that
\[
  \nu(F_1)>\tfrac12\nu(F_0),
  \qquad
  \lim_{n\to\infty}\sup_{z\in F_1}
  \|\widetilde{\mathcal C}_{1/n}\mu(z)-D(z)\|_2=0.
\]
By inner regularity, choose a compact set \(F\subset F_1\) with
\(\nu(F)>\tfrac12\nu(F_1)>0\). Then
\[
  0\le \sup_{z\in F}
  \|\widetilde{\mathcal C}_{1/n}\mu(z)-D(z)\|_2
  \le \sup_{z\in F_1}
  \|\widetilde{\mathcal C}_{1/n}\mu(z)-D(z)\|_2
  \longrightarrow 0.
\]
Since \(F\subset F_0\subset F_{i,k}\), we have \(|\mu_s|(F)=0\) and
\[
  k^{-1}\le \operatorname{Re}\langle f(z),v_i\rangle
  \le \|f(z)\|\le k,
  \qquad z\in F.
\]
Restricting \(f\) to \(F\) and taking \(v=v_i\) and \(d=k^{-1}\)
completes the proof.
\end{proof}

\begin{proposition}\label{gen:local}
Let \(U\subset\mathbb C\) be a bounded open disk with
\(\gamma(U\cap\Lambda)>0\). Then there exist
\(x\in H\setminus\{0\}\) and a constant \(C_x<\infty\) such that
\[
\|p(T)x\|\le C_x\sup_{z\in U}|p(z)|,
\qquad p\in\mathbb C[z].
\]
\end{proposition}

\begin{proof}
Since
\[
0<\gamma(U\cap\Lambda)
=\sup\{\gamma(L):L\subset U\cap\Lambda,\ L\text{ compact}\},
\]
there is a compact set \(L_0\subset U\cap\Lambda\) with
\(\gamma(L_0)>0\). By Lemma~\ref{cau:capacity}\textup{(2)}, there is
a finite positive measure \(\nu\) supported on \(L_0\) satisfying
\eqref{cau:measure}. In particular, \(\nu(L_0)\ge c\gamma(L_0)>0\).
For this choice of \(L_0\) and \(\nu\), use the notation
\(\mathcal M,Z,\mathcal S,\mathcal V\) introduced before
Lemma~\ref{gen:compact-selection}.

Suppose first that \(\nu(\mathcal S\cup\mathcal V)>0\). Choose
\(F,e,v,d,f\) as in Lemma~\ref{gen:compact-selection}. Since
\(\nu(F)>0\) and \(f\) is bounded and Borel on \(F\),
Lemma~\ref{gen:cutoff} gives a measurable function
\(h:F\to[0,1]\) such that
\[
\int_F h\,d\nu\ge\tfrac12\nu(F),\qquad
\sup_{\substack{w\in\mathbb C\\ r>0}}
\left\|\int_F k_r(z-w)h(z)f(z)\,d\nu(z)\right\|<\infty.
\]
Choose an \([0,1]\)-valued Borel representative of \(h\) on \(F\),
and extend it by zero to \(\mathbb C\setminus F\). \textcolor{black}{Set} \(g=hf\) on \(F\) and \(u=h(N)Ae\). Then
\[
E(F)u=(\mathbf1_F h)(N)Ae=u,\qquad
\sup_{n\ge1,\,w\in\mathbb C}
\left\|\int_F k_{1/n}(z-w)g(z)\,d\nu(z)\right\|<\infty.
\]
Since \(|\mu_s|(F)=0\) and \(\Delta f=\mathcal M e\) on \(F\),
\eqref{model:borel-calculus} gives, for every \(p\in\mathbb C[z]\),
\[
\begin{aligned}
Bp(N)u
&=B(ph)(N)Ae=\left(\int_F ph\,d\mu\right)e\\
&=\int_F ph\,\mathcal M e\,d\nu
=\int_F p\,\Delta g\,d\nu.
\end{aligned}
\]
Thus \eqref{gen:correction} holds. By Lemma~\ref{gen:compact-selection},
\(\widetilde{\mathcal C}_{1/n}\mu\) converges uniformly to \(D\) on \(F\)
in the Hilbert--Schmidt norm. Applying
Lemma~\ref{gen:polynomial-moments} with \(r_n=1/n\), we obtain \(x\in H\)
and \(C_x<\infty\) such that, for every \(p\in\mathbb C[z]\),
\[
\|p(T)x\|\le C_x\sup_{z\in U}|p(z)|,\qquad
Bp(T)x=\int_F p\,g\,d\nu.
\]
Taking \(p=1\) gives
\[
\operatorname{Re}\langle Bx,v\rangle
=\int_F h\,\operatorname{Re}\langle f,v\rangle\,d\nu
\ge d\int_F h\,d\nu\ge\tfrac12 d\,\nu(F)>0.
\]
Hence \(x\ne0\).

Suppose now that \(\nu(\mathcal S\cup\mathcal V)=0\). Since
\(\nu(Z)=0\), we have \(\nu(Z\cup\mathcal S\cup\mathcal V)=0\).
For every \(z\notin Z\cup\mathcal S\cup\mathcal V\), the definitions
of \(\mathcal S\) and \(\mathcal V\) give that \(\Delta(z)\) is
invertible and \(\mathcal M(z)=0\). Hence
\[
d\mu=\mathcal M\,d\nu+d\mu_s=d\mu_s,
\qquad |\mu|=|\mu_s|\perp\nu.
\]
If \(a\in\operatorname{supp}\nu\), then every ball \(B(a,t)\),
\(t>0\), has positive \(\nu\)-measure. Together with
\eqref{cau:measure}, this gives
\[
0<\nu(B(a,t))\le t \qquad (t>0).
\] Since \(|\mu|\perp\nu\), the differentiation theorem for finite Borel
measures on \(\mathbb C\) gives
\[
\lim_{t\to0^+}\frac{|\mu|(B(a,t))}{\nu(B(a,t))}=0
\qquad\text{for \(\nu\)-almost every }a\in\operatorname{supp}\nu.
\]
For each such \(a\), the bound \(0<\nu(B(a,t))\le t\) yields
\[
0\le\frac{|\mu|(B(a,t))}{t}
\le\frac{|\mu|(B(a,t))}{\nu(B(a,t))}
\longrightarrow0,
\qquad t\to0^+.
\]
Since \(\nu\) is supported on \(L_0\) and \(\nu(Z\cup\mathcal S)=0\),
\[
\nu\bigl(\operatorname{supp}\nu\setminus(Z\cup\mathcal S)\bigr)
=\nu(L_0)>0.
\]
By the preceding limit, we may therefore choose
\(a\in\operatorname{supp}\nu\setminus(Z\cup\mathcal S)\) such that
\[
\lim_{t\to0^+}\frac{|\mu|(B(a,t))}{t}=0.
\]
Then \(a\in L_0\subset U\cap\Lambda\), and \(\Delta(a)\) is invertible
by the definition of \(\mathcal S\). Moreover, since \(a\notin Z\),
\[
\|\mathcal C_t\mu(a)-D(a)\|_2\longrightarrow0
\qquad (t\to0^+).
\]
Lemma~\ref{model:regular}, applied at \(a\), now gives the assertion.
\end{proof}

\begin{lemma}\label{model:subnormal-singular}
With the notation of Lemma~\ref{model:subnormal-extension}, let
\(U\subset\mathbb C\) be a bounded open disk and set
\[
\mathcal V_U
=\{z\in U:I_\mu(z)<\infty \text{ and } \Delta(z)\text{ is singular}\}.
\]
If \(m_2(\mathcal V_U)>0\), there exist \(v\in H\setminus\{0\}\) and
\(C_v<\infty\) such that
\[
\|p(T^*)v\|\le C_v\sup_{z\in U}|p(z)|,
\qquad p\in\mathbb C[z].
\]
\end{lemma}

\begin{proof}
Since \(I_\mu\) is lower semicontinuous, \(\Delta\) is norm-Borel,
and the singular operators form a norm-closed subset of
\(\mathcal B(\mathcal E)\), the set \(\mathcal V_U\) is Borel. For each \(z\in\mathcal V_U\), the operator \(\Delta(z)\) is singular by
definition and Fredholm of index zero by
Lemma~\ref{model:regularized-transform}\textup{(1)}. Hence \(\ker\Delta(z)\ne\{0\}\). Let \(Q(z)\) be the
orthogonal projection onto \(\ker\Delta(z)\). For each \(b\in\mathcal E\),
the spectral theorem gives
\[
Q(z)b=\lim_{k\to\infty}
(I_{\mathcal E}+k\Delta(z)^*\Delta(z))^{-1}b,
\qquad z\in\mathcal V_U.
\]
Since \(\Delta\) is norm-Borel, each map on the right is Borel in \(z\).
Thus \(z\mapsto Q(z)b\) is Borel on \(\mathcal V_U\).

Let \(\{b_j\}_{j=1}^{\infty}\) be a countable orthonormal basis of \(\mathcal E\).
For each \(z\in\mathcal V_U\), \(Q(z)\ne0\), so \(Q(z)b_j\ne0\) for some \(j\).
Consequently,
\[
\mathcal V_U=\bigcup_j\{z\in\mathcal V_U:Q(z)b_j\ne0\}.
\]
Since \(m_2(\mathcal V_U)>0\), countable subadditivity gives
\(m_2\{z\in\mathcal V_U:Q(z)b_{j_0}\ne0\}>0\) for some \(j_0\).
Fix this index and \textcolor{black}{let} \(e=b_{j_0}\). Then
\[
\int_{\mathcal V_U}\|Q(z)e\|^2\,dm_2(z)>0.
\]
\textcolor{black}{Set} \(g(z)=Q(z)e\) on \(\mathcal V_U\) and \(g(z)=0\) on
\(U\setminus\mathcal V_U\).
For \(z\in\mathcal V_U\), the orthogonal projection \(Q(z)\) onto
\(\ker\Delta(z)\) satisfies \(\|Q(z)e\|\le\|e\|=1\) and
\(\Delta(z)Q(z)e=0\). Since \(g=0\) on \(U\setminus\mathcal V_U\), we obtain
\(\|g\|_\infty\le1\) and \(\Delta g=0\) on \(U\). For each
\(p\in\mathbb C[z]\), apply \eqref{sub:transfer} to \(pg\) to obtain
\[
WJ(pg)
=J(zpg)-A\int_U p(z)\Delta(z)g(z)\,dm_2(z)
=J(zpg).
\]
By induction, \(W^kJg=J(z^kg)\) for every \(k\ge0\). By linearity,
\(p(W)Jg=J(pg)\) for every \(p\in\mathbb C[z]\).

\textcolor{black}{Let} \(v=PJg\). Since \(B=B_0P\), we have \(B_0v=BJg\).
Also, \(D=I_{\mathcal E}-\Delta\) gives \(Dg=g\). Hence
\eqref{sub:transfer-bound} gives
\[
\langle B_0v,e\rangle
=\langle BJg,e\rangle
=\int_{\mathcal V_U}\langle Q(z)e,e\rangle\,dm_2(z)
=\int_{\mathcal V_U}\|Q(z)e\|^2\,dm_2(z)>0.
\]
Thus \(v\ne0\). Finally, by \eqref{sub:extension} and
\eqref{sub:transfer-bound},
\[
\|p(T^*)v\|
=\|Pp(W)Jg\|
=\|PJ(pg)\|
\le c_U\|A\|_2\|g\|_\infty\sup_{z\in U}|p(z)|,
\quad p\in\mathbb C[z].
\]
\end{proof}

\begin{proposition}\label{sub:local}
Let \(S\) be a bounded subnormal operator on \(H\), let \(K\in S_1(H)\),
and \textcolor{black}{set} \(T=S+K\). Let \(\sigma_{\mathrm e}(S)\) be the Fredholm essential spectrum
of \(S\), and \textcolor{black}{set} \(\sigma_{\mathrm e}(S)^\#=\{\bar\lambda:\lambda\in\sigma_{\mathrm e}(S)\}\).
For a bounded open disk \(U\subset\mathbb C\), suppose that
\begin{enumerate}
\item[\textup{(1)}] \(\ker(T-\lambda I)=\ker(T^*-\bar\lambda I)=\{0\}\)
for every \(\lambda\in\mathbb C\).
\item[\textup{(2)}] \(m_2(U\cap\sigma_{\mathrm e}(S)^\#)>0\).
\end{enumerate}
Then there exist \(v\in H\setminus\{0\}\) and \(C_v<\infty\) such that
\begin{equation}\label{sub:local-bound}
\|p(T^*)v\|\le C_v\sup_{z\in U}|p(z)|,
\qquad p\in\mathbb C[z].
\end{equation}
\end{proposition}

\begin{proof}
Let \(\mathcal V_U\) be the set defined in Lemma~\ref{model:subnormal-singular}
for the disk \(U\). If \(m_2(\mathcal V_U)>0\), that lemma gives
\eqref{sub:local-bound}.
Suppose that \(m_2(\mathcal V_U)=0\). By Lemma~\ref{cau:area-lemma},
\(I_\mu\) is finite \(m_2\)-almost everywhere on \(U\). Hence, by
hypothesis \textup{(2)},
\[
m_2\bigl(\{z\in(U\cap\sigma_{\mathrm e}(S)^\#)\setminus\mathcal V_U:
I_\mu(z)<\infty\}\bigr)
=m_2(U\cap\sigma_{\mathrm e}(S)^\#)>0.
\]
Choose \(a\) in this set. Since \(I_\mu(a)<\infty\) and \(a\notin\mathcal V_U\),
\(\Delta(a)\) is invertible.

\textcolor{black}{Let} \(\lambda=\bar a\).
Then \(\lambda\in\sigma_{\mathrm e}(S)\), so \(S-\lambda I\) is not Fredholm.
Hence
\(T-\lambda I=(S-\lambda I)+K\) is not Fredholm.
By \textup{(1)}, \(T-\lambda I\) is injective with dense range.
Since \(T-\lambda I\) is not Fredholm, its range cannot be closed.

Let \(M\) be a finite-dimensional subspace of \(H\). If \(T-\lambda I\)
were bounded below on \(M^\perp\), then \((T-\lambda I)M^\perp\)
would be closed. Since \((T-\lambda I)M\) is finite-dimensional,
\[
\operatorname{ran}(T-\lambda I)
=(T-\lambda I)M^\perp+(T-\lambda I)M
\]
would also be closed, a contradiction. Hence
\[
\inf_{\substack{s\in M^\perp\\ \|s\|=1}}\|(T-\lambda I)s\|=0.
\]

Taking \(M=\{0\}\), choose a unit vector \(s_1\in H\) with
\(\|(T-\lambda I)s_1\|<1\). Once \(s_1,\ldots,s_{n-1}\) have been
chosen, take \(M=\operatorname{span}\{s_1,\ldots,s_{n-1}\}\) and
choose a unit vector \(s_n\in M^\perp\) with
\(\|(T-\lambda I)s_n\|<1/n\). The resulting sequence \((s_n)\) is
orthonormal, so \(s_n\rightharpoonup0\). Since \(K\) is compact,
\(\|Ks_n\|\to0\). Thus \(S=T-K\) and \(\lambda=\bar a\) give
\begin{equation}\label{sub:singular-sequence}
\|(S-\bar a I)s_n\|
\le \|(T-\lambda I)s_n\|+\|Ks_n\|
\le \frac{1}{n}+\|Ks_n\|\longrightarrow0.
\end{equation}
Since \(I_\mu(a)<\infty\), \(\Delta(a)\) is invertible, and
\eqref{sub:singular-sequence} holds, Lemma~\ref{model:subnormal-regular}
applied at \(a\) gives \eqref{sub:local-bound}.
\end{proof}

\section{Proofs of the main theorems}\label{sec:half}

This section gives the proofs of the main theorems. For any compact set $\mathcal K\subset \mathbb C$ and  any polynomial $p\in\mathbb C[z]$, let 
$$
\|p\|_{\mathcal K}=\max\{|p(z)|; z\in \mathcal K\}.
$$
Let $D$ be a closed disk of positive radius. For any $R\in \mathcal B(H)$, set
$$
  c_D(x;R)=
  \sup\{\|p(R)x\|;\,p\in\mathbb C[z],\ \|p\|_D\le1\},
 $$
 and let
 \[
  X_R(D)=\{x\in H:c_D(x;R)<\infty\}.
\]

For \(F\subset\mathbb C\), write \(F^\#=\{\bar z:z\in F\}\); 
moreover, for any polynomial \(p(z)=\sum\limits_{j=0}^n a_jz^j\), set
$$p^\#(z)=\overline{p(\bar z)}=\sum_{j=0}^n\overline{a_j}z^j.$$
Then \(\|p^\#\|_{D^\#}=\|p\|_D\) for every closed disk \(D\), and
\[
  p(R)^*
  =\left(\sum_{j=0}^n a_jR^j\right)^*
  =\sum_{j=0}^n\overline{a_j}(R^*)^j
  =p^\#(R^*).
\]
The following lemma is inspired by \cite[Section 1.2]{GallardoGonzalez2022}.

\begin{lemma}\label{eng:separate}
Let \(D_1,D_2\subset\mathbb C\) be disjoint closed disks. Then for any \(R\in\mathcal B(H)\),
\[
  X_R(D_1)\perp X_{R^*}(D_2^\#).
\]
Moreover, if \(D_1,\ldots,D_m\) are pairwise disjoint closed disks
and the vectors \(x_j\in X_R(D_j)\) satisfy \(\sum\limits_{j=1}^m x_j=0\),
then \(x_1=\cdots=x_m=0\).
\end{lemma}

\begin{proof}
The function \(\mathbf1_{D_1}\) on \(D_1\cup D_2\) extends holomorphically to a neighborhood of this union. Since the union has connected complement, Runge's theorem gives polynomials \(q_n\) such that
\[
 \|q_n-\mathbf1_{D_1}\|_{D_1\cup D_2}
 =\max\{\|1-q_n\|_{D_1},\|q_n\|_{D_2}\}
 \xrightarrow{n\to\infty}0.
\]
Let \(x\in X_R(D_1)\) and \(y\in X_{R^*}(D_2^\#)\). Using
\(q_n(R)^*=q_n^\#(R^*)\), we obtain
\[
  \langle x,y\rangle
  =\langle (I-q_n(R))x,y\rangle+\langle x,q_n^\#(R^*)y\rangle.
\]
Since \(\|q_n^\#\|_{D_2^\#}=\|q_n\|_{D_2}\), it follows that
\[
\begin{aligned}
  |\langle x,y\rangle|
  &\le\lim\limits_{n\to\infty}\left( \|(I-q_n(R))x\|\,\|y\|+\|x\|\,\|q_n^\#(R^*)y\|\right)\\
  &\le\lim\limits_{n\to\infty}\left( c_{D_1}(x;R)\,\|1-q_n\|_{D_1}\,\|y\|
      +\|x\|\,c_{D_2^\#}(y;R^*)\,\|q_n\|_{D_2}\right)=0.
\end{aligned}
\]
It follows that \(X_R(D_1)\perp X_{R^*}(D_2^\#)\).

Next, by the hypothesis,  $\{D_j\}_{j=1}^m$ are disjoint,    
and \(x_j\in X_R(D_j)\) satisfy \(\sum\limits_{j=1}^m x_j=0\). Fix
\(i\in\{1,\ldots,m\}\). By Runge's theorem again, there are polynomials
\(p_n\) such that
\[
 \lim\limits_{n\to\infty}\max_{1\le j\le m}
  \|p_n-\delta_{ij}\|_{D_j}=0.
\]
Since \(p_n(R)\sum\limits_{j=1}^m x_j=0\), we have
\(x_i=\sum\limits_{j=1}^m(\delta_{ij}I-p_n(R))x_j\). For simplicity, write $\varepsilon_n=\max\limits_{1\le j\le m}
  \|p_n-\delta_{ij}\|_{D_j}$, and we have
\[
  \|x_i\|
  \le\lim\limits_{n\to\infty} \sum_{j=1}^m c_{D_j}(x_j;R)\,\|\delta_{ij}-p_n\|_{D_j}
  \le\lim\limits_{n\to\infty} \varepsilon_n\sum_{j=1}^m c_{D_j}(x_j;R)
 = 0.
\]
This completes the proof. 
\end{proof}

\begin{proposition}\label{eng:positive}
Let \(N\) be normal and \(K\in S_1(H)\).
If \(\gamma(\sigma(N))>0\), then \(T=N+K\) has a hyperinvariant
half-space \(\mathcal L\subset H\).
\end{proposition}

\begin{proof}
By Lemma~\ref{cau:capacity}\textup{(2)}, there is a nonzero finite positive measure
\(\nu\) supported on \(\Lambda=\sigma(N)\) and a constant $C>0$ such that
\(\nu(F)\le C\gamma(F)\) for every compact \(F\subset\mathbb C\).
For each \(a\in\mathbb C\),
\[
  0\le\nu(\{a\})\le C\gamma(\{a\})=0.
\]
Thus \(\nu\) has no atoms. If \(a\in\operatorname{supp}\nu\), then  for every \(r>0\),
$$\nu(B(a,r)\setminus\{a\})=\nu(B(a,r))>0.$$
Therefore \(\operatorname{supp}\nu\) has no isolated points. Let \(U\) be an open disk which intersects with \(\operatorname{supp}\nu\).
By inner regularity of the measure $\nu$,
\[
  0<\nu(U)
  =\sup\{\nu(F); F\subset U\cap\Lambda \text{ and } F\ \mathrm{compact}\}.
\]
Choose an \(F\) with \(\nu(F)>0\). Then
\(\gamma(F)\ge C^{-1}\nu(F)>0\). By Proposition~\ref{gen:local}
there is an \(x\ne0\) and  \(C_x<\infty\) such that for any $p\in\mathbb C[z]$
\[
  \|p(T)x\|\le 
  C_x\|p\|_{\overline U}.
\]
Consequently \(c_{\overline U}(x;T)\le C_x\), and hence
\(X_T(\overline U)\ne\{0\}\). Reflection gives that for $ F^\#\subset U^\#\cap\Lambda^\#
         =U^\#\cap\sigma(N^*),$
\[
  \gamma(F^\#)=\gamma(F)>0.
\]
Since \(N^*\) is normal and \(K^*\in S_1(H)\),
applying Proposition~\ref{gen:local} to \(T^*\), we have
\(y\ne0\) and \(C_y<\infty\) such that for any $p\in\mathbb C[z]$,
\[
  \|p(T^*)y\|\le C_y\sup_{U^\#}|p|
  =C_y\|p\|_{\overline U^\#}.
\]
Thus \(c_{\overline U^\#}(y;T^*)\le C_y\), and
\(X_{T^*}(\overline U^\#)\ne\{0\}\).

Choose distinct \(a_1,a_2\in\operatorname{supp}\nu\), and set $r=|a_1-a_2|/3$ 
 and let \(L_i=\overline{B(a_i,r)}\). Then 
\(L_1\cap L_2=\varnothing\). Since \(\operatorname{supp}\nu\) has
no isolated points, for each \(i=1,2\) there are
\(a_{ij}\in\operatorname{supp}\nu\setminus\{a_i\}\) such that
\[
  d_{i1}<r/2 \text{ and } 0<d_{i,j+1}<d_{ij}/4,
\]
where $ d_{ij}:=|a_{ij}-a_i|$.
Let \(D_{ij}=\overline{B(a_{ij},d_{ij}/4)}\). For \(z\in D_{ij}\),
\[
  |z-a_i|\le |z-a_{ij}|+d_{ij}
  \le \tfrac54d_{ij}\le\tfrac54d_{i1}<\tfrac58r.
\]
Therefore \(D_{ij}\subset B(a_i,r)\). For \(k>j\) we have
\(d_{ik}<d_{ij}/4\), whence
\[
\begin{aligned}
  \operatorname{dist}(D_{ij},D_{ik})
  \ge |a_{ij}-a_{ik}|-\tfrac14(d_{ij}+d_{ik})
  \ge \tfrac34d_{ij}-\tfrac54d_{ik}
   >\tfrac7{16}d_{ij}>0.
\end{aligned}
\]
Thus, for each \(i\), the disks \(D_{ij}\) are pairwise disjoint. For each \(j\ge1\), choose
\[
  0\ne x_j\in X_T(D_{1j})\subset X_T(L_1), \text{ and }
  0\ne y_j\in X_{T^*}(D_{2j}^\#)\subset X_{T^*}(L_2^\#).
\]
By Lemma~\ref{eng:separate}, each family is linearly independent and
\(X_T(L_1)\perp X_{T^*}(L_2^\#)\). \textcolor{black}{Let}
\(\mathcal L=\overline{X_T(L_1)}\). Then
\[
  \operatorname{span}\{x_j:j\ge1\}\subset\mathcal L,\text{ and }
  \operatorname{span}\{y_j:j\ge1\}\subset\mathcal L^\perp.
\]
It follows that \(\dim\mathcal L=\dim\mathcal L^\perp=\infty\).

At last, for \(S\in\{T\}'\) and \(x\in X_T(L_1)\), 
\[
  c_{L_1}(Sx;T)
  =\sup
     \{\|Sp(T)x\|;p\in\mathbb C[z],\ \|p\|_{L_1}\le1\}
  \le \|S\|\,c_{L_1}(x;T)<\infty.
\]
Thus \(SX_T(L_1)\subset X_T(L_1)\), and continuity yields
\[
  S\mathcal L
  \subset\overline{SX_T(L_1)}
  \subset\overline{X_T(L_1)}
  =\mathcal L
  \qquad(S\in\{T\}').
\]
Hence \(\mathcal L\) is a hyperinvariant half-space for \(T\).
\end{proof}

We now prove the main theorems.

\begin{proof}[Proof of Theorem~\ref{main:invariant}]
\textcolor{black}{Let} \(T=S+K\). If \(S\) is nonnormal, Theorem~\ref{sub:nonnormal-intransitive}
below gives the conclusion. Suppose that \(S\) is normal, and write
\(N=S\) and \(\Lambda=\sigma(N)\). Since \(K\) is compact,
\[
  \sigma_{\mathrm e}(T)=\sigma_{\mathrm e}(N),\qquad
  \operatorname{ind}(T-zI)=\operatorname{ind}(N-zI)=0
  \quad(z\notin\sigma_{\mathrm e}(N)).
\]
If \(z\in\sigma(T)\setminus\sigma_{\mathrm e}(N)\), then \(T-zI\) is Fredholm
of index zero and is not invertible. Hence
\[
  \dim\ker(T-zI)=\dim\ker((T-zI)^*)>0.
\]
An eigenvector spans an invariant line, so \(T\) is intransitive.
We may assume that \(\sigma(T)=\sigma_{\mathrm e}(N)\).

If \(\sigma_{\mathrm e}(N)\) is disconnected, a Riesz projection associated with
either part of a separation of \(\sigma_{\mathrm e}(N)\) has nonzero proper
invariant range, so \(T\) is intransitive.

If \(\sigma_{\mathrm e}(N)=\{\lambda\}\), normality of \(N\) implies that
\(N-\lambda I\) is compact. Thus \(T-\lambda I\) is compact.
If \(T=\lambda I\), every line is invariant. Otherwise
\(T-\lambda I\) is intransitive by the Aronszajn--Smith theorem
\cite{AS}. In either case \(T\) is intransitive.

In the remaining case, \(\sigma_{\mathrm e}(N)\) is connected and
\(\operatorname{diam}(\sigma_{\mathrm e}(N))>0\). We claim that
\[
  \gamma(\sigma_{\mathrm e}(N))\ge\tfrac14\operatorname{diam}(\sigma_{\mathrm e}(N)).
\]
Let \(\Omega\) be the component of
\(\widehat{\mathbb C}\setminus\sigma_{\mathrm e}(N)\) containing \(\infty\).
By the Riemann mapping theorem, there is a conformal map
\[
  g:\{z\in\mathbb C:|z|>1\}\cup\{\infty\}\longrightarrow\Omega,
  \text{ and } g(\infty)=\infty,
\]
normalized by \(g(z)=cz+c_0+O(z^{-1})\) as \(z\to\infty\), with
\(c>0\). Define \(f(z)=1/g^{-1}(z)\) on \(\Omega\) and \(f=0\) on the other
components of \(\widehat{\mathbb C}\setminus\sigma_{\mathrm e}(N)\). Therefore
\[
  f\in H^\infty(\widehat{\mathbb C}\setminus\sigma_{\mathrm e}(N)),\, 
  \|f\|_\infty\le1,\, 
  f(\infty)=0,\text{ and } f'(\infty)=c.
\]
Thus \(\gamma(\sigma_{\mathrm e}(N))\ge c\). For distinct \(a,b\in\sigma_{\mathrm e}(N)\), set
\[
  F_a(\zeta)=\frac{c}{g(1/\zeta)-a}\quad(0<|\zeta|<1),
  \text{ and } F_a(0)=0.
\]
The function \(F_a\) is univalent in \(\mathbb D\), and the expansion
of \(g\) gives \(F_a(\zeta)=\zeta+O(\zeta^2)\). Since \(b\notin\Omega\),
the Koebe quarter theorem yields
\[
  \frac{c}{b-a}\notin F_a(\mathbb D),\qquad
  B(0,\tfrac14)\subset F_a(\mathbb D),
  \qquad |b-a|\le4c.
\]
Taking the supremum over \(a,b\in\sigma_{\mathrm e}(N)\), we obtain
\[
  \gamma(\Lambda)\ge\gamma(\sigma_{\mathrm e}(N))\ge c
  \ge\tfrac14\operatorname{diam}(\sigma_{\mathrm e}(N))>0.
\]
By Proposition~\ref{eng:positive}, \(T\) has a hyperinvariant
half-space and hence is intransitive.
\end{proof}
For a compact set \(Y\subset\mathbb C\), let \(Q(Y)\) be the rational
functions with poles off \(Y\), and let \(R(Y)\) be their uniform closure
in \(C(Y)\).

\begin{theorem}\label{sub:nonnormal-intransitive}
Let \(S\) be a nonnormal subnormal operator on \(H\) and
\(K\in S_1(H)\). Then \(T=S+K\) is intransitive.
\end{theorem}

\begin{proof}
If \(T\) has an eigenvector, its span is invariant. If \(T^*\) has an
eigenvector, the orthogonal complement of its span is invariant for \(T\).
In either case \(T\) is intransitive. We may therefore assume that
\(\ker(T-\lambda I)=\ker(T^*-\bar\lambda I)=\{0\}\) for all
\(\lambda\in\mathbb C\).

\textcolor{black}{Let} \(X=\sigma(S)\). Since \(K\) is compact,
\(\sigma_{\mathrm e}(S)=\sigma_{\mathrm e}(T)\). For
\(\lambda\notin\sigma_{\mathrm e}(S)\), the Fredholm operator \(T-\lambda I\)
has trivial kernel and dense range, hence is invertible. It follows that
\[
  \operatorname{ind}(S-\lambda I)=\operatorname{ind}(T-\lambda I)=0.
\]
For \(\lambda\in X\setminus\sigma_{\mathrm e}(S)\), the operator
\(S-\lambda I\) is thus Fredholm of index zero and singular. Hence
\[
  \dim\ker(S-\lambda I)=\dim\ker(S^*-\bar\lambda I)>0.
\]
Since \(S\) is subnormal, \(S-\lambda I\) is hyponormal and
\(\ker(S-\lambda I)\subseteq\ker(S^*-\bar\lambda I)\). Consequently,
\[
  \ker(S-\lambda_1 I)\perp\ker(S-\lambda_2 I)
  \qquad(\lambda_1\ne\lambda_2).
\]
Since \(H\) is separable, \(X\setminus\sigma_{\mathrm e}(S)\) is countable.

Next, we claim that $R(X)\ne C(X)$.
In fact, let $N$ be the minimal normal extension of $S$, with spectral measure $E_N$. For $\lambda\notin X$, choose $r>0$ such that $r\|(S-\lambda I)^{-1}\|<1$.
For $x\in H$ and $m\ge1$,
\[
\begin{aligned}
 \|E_N(B(\lambda,r))x\|
 &=\|E_N(B(\lambda,r))(N-\lambda I)^m(S-\lambda I)^{-m}x\|\\
 &\le \bigl(r\|(S-\lambda I)^{-1}\|\bigr)^m\|x\|.
\end{aligned}
\]
Letting \(m\to\infty\) gives \(E_N(B(\lambda,r))|_H=0\). Since
\[
  E_N(B(\lambda,r))N^*=N^*E_N(B(\lambda,r)),
\]
minimality of the normal extension of \(S\) gives \(E_N(B(\lambda,r))=0\).
Thus \(\sigma(N)\subseteq X\), and \(q(N)|_H=q(S)\) for every \(q\in Q(X)\).
If \(R(X)=C(X)\), there are \(q_n\in Q(X)\) converging uniformly on \(X\)
to \(z\mapsto\bar z\). Functional calculus for \(N\) gives
\(q_n(N)\to N^*\) in norm. Since \(q_n(N)H\subseteq H\), we have
\(N^*H\subseteq H\). Thus \(H\) reduces \(N\) and \(S\) is normal,
a contradiction.

Let \(\varrho=\mathbf1_X\,dm_2\) and \(F=\operatorname{supp}\varrho\).
Since \(\sigma_{\mathrm e}(S)\) is closed and
\(X\setminus\sigma_{\mathrm e}(S)\) is countable,
\[
  \varrho(\mathbb C\setminus\sigma_{\mathrm e}(S))=0,
  \qquad F\subseteq\sigma_{\mathrm e}(S).
\]
Since \(R(X)\ne C(X)\), the Hahn--Banach and Riesz representation theorems
give a nonzero finite complex measure \(\vartheta\) on \(X\) such that
\[
  \int_X r(z)\,d\vartheta(z)=0
  \qquad(r\in R(X)).
\]
Its Cauchy transform \(\mathcal C\vartheta\in L^1_{\mathrm{loc}}(\mathbb C)\)
satisfies
\[
  \mathcal C\vartheta(z)=\int_X\frac{d\vartheta(w)}{z-w}=0
  \qquad(z\notin X).
\]
Since \(m_2(X\setminus F)=0\), we have \(\mathcal C\vartheta=0\) almost
everywhere on \(\mathbb C\setminus F\). Thus, in the sense of distributions,
\[
  \bar\partial\mathcal C\vartheta=\pi\vartheta,
  \qquad \operatorname{supp}\vartheta\subseteq F.
\]
Consequently, \(\mathcal C\vartheta\) is holomorphic and identically zero
on \(\mathbb C\setminus F\). Differentiating at \(a\notin F\) gives
\[
  \int_F\frac{d\vartheta(w)}{(w-a)^j}
  =-\frac{(\mathcal C\vartheta)^{(j-1)}(a)}{(j-1)!}=0
  \qquad(j\ge1).
\]
Also, \(\int_F w^n\,d\vartheta(w)=0\) for every \(n\ge0\).
By partial fractions and uniform approximation, \(\vartheta\) annihilates
\(R(F)\). Since \(\vartheta\ne0\), we have \(R(F)\ne C(F)\). {\color{black} By the thickness criterion in~\cite{Brown1987}, as stated in~\cite[p.~361]{Lauric2001}}, there exists a nonempty bounded open set
\(\Omega\subset\mathbb{C}\) such that
\[
  \sup_{z\in\Omega}|f(z)|
  =\sup_{z\in\Omega\cap F}|f(z)|
  \qquad (f\in H^\infty(\Omega)).
\]
Thus \(F\) is thick.
For \(g\in H^\infty(\Omega^\#)\), apply the preceding identity to
\(f(z)=\overline{g(\bar z)}\). Then
\[
  \sup_{z\in\Omega^\#}|g(z)|
  =\sup_{z\in\Omega^\#\cap F^\#}|g(z)|,
\]
and hence \(F^\#\) is thick as well.

Fix \(c\in F^\#\), and let \(\mathcal O\) be an open neighborhood of \(c\).
Choose \(t>0\) such that \(U=B(c,t)\) satisfies \(\overline U\subset\mathcal O\).
By the definition of \(F\) and the countability of \(X\setminus\sigma_{\mathrm e}(S)\),
\[
  0<m_2(U\cap X^\#)
   =m_2(U\cap X^\#\cap\sigma_{\mathrm e}(S)^\#)
   \le m_2(U\cap\sigma_{\mathrm e}(S)^\#).
\]
Proposition~\ref{sub:local} yields \(v\in H\setminus\{0\}\) and \(C_v>0\) such that for any $p\in\mathbb C[z]$
\[
  \|p(T^*)v\|\le C_v\sup_{z\in U}|p(z)|.
\]
By taking \(p(z)=(z-c)^j\), we have  for every \(j\ge0\)
\(\|(T^*-cI)^jv\|\le C_v t^j\).
Hence for $|\lambda-c|>t$, the series
\begin{equation}\label{sub:local-resolvent}
  q(\lambda)=\sum_{j=0}^{\infty}
    \frac{(T^*-cI)^jv}{(\lambda-c)^{j+1}}
\end{equation}
converges locally uniformly on \(\mathbb C\setminus\overline U\). Therefore, for $\lambda\in\mathbb C\setminus\overline U$,
by multiplying  \(\lambda I-T^*\) on both sides of \eqref{sub:local-resolvent}, we have
\[
  (\lambda I-T^*)q(\lambda)=v.
\]
Thus the local spectrum of \(T^*\) at \(v\) satisfies
\(\sigma_{T^*}(v)\subseteq\overline U\subset\mathcal O\).

The series \eqref{sub:local-resolvent} gives a holomorphic solution on
\(\mathbb C\setminus\overline U\). Since \(\mathcal O\) was arbitrary,
\(c\) belongs to the localizable spectrum of \(T^*\). This spectrum contains the thick set \(F^\#\) and is therefore
thick. By
\cite[Corollary~3.4]{EschmeierPrunaru2002}, there is a closed subspace
\(\{0\}\subsetneq M\subsetneq H\) with \(T^*M\subseteq M\).
Then \(M^\perp\) is nonzero, proper and invariant for \(T\), so
\(T\) is intransitive.
\end{proof}
Recall that a compact set has positive Riesz \(1\)-capacity if it
supports a probability measure \(\omega\) with \(I_1(\omega)<\infty\).

\begin{theorem}\label{main:half1}
Let $N$ be normal and $K\in S_1(H)$. If $\sigma(N)$
has positive Riesz $1$-capacity, then $N+K$ has a hyperinvariant
half-space.
\end{theorem}

\begin{proof}
\textcolor{black}{Set} \(T=N+K\) and \(\Lambda=\sigma(N)\), and use the normal-operator
notation of Section~\ref{sec:model}. Choose a probability measure
\(\omega\) supported on \(\Lambda\)
with \(I_1(\omega)<\infty\). By Lemma~\ref{cau:energy}, there is
a compact set \(F\subset\Lambda\) with \(\omega(F)>0\) such that, for \(\tau=\omega|_F\),
\[
  b_\tau(r)\xrightarrow{r\to0^+}0,\qquad C_\tau<\infty.
\]
Fix a bounded open disk \(U\) with \(\tau(U)>0\), and \textcolor{black}{let}
\(\eta=\tau|_U\). By Tonelli's theorem,
\[
\begin{aligned}
  \int_U I_\mu(z)\,d\eta(z)
  =\int_{\mathbb C}
    \left(\int_U\frac{d\tau(z)}{|z-w|}\right)d|\mu|(w)
  \le C_\tau|\mu|(\mathbb C)<\infty.
\end{aligned}
\]
Thus \(I_\mu(a)<\infty\) for \(\eta\)-almost every \(a\in U\cap\Lambda\).
At each such point, dominated convergence gives
\[
  \lim\limits_{\varepsilon\to0^+}\mathcal C_\varepsilon\mu(a)=
  \int_{\mathbb C}\frac{d\mu(w)}{a-w}=D(a)
  \quad\text{in }S_2(\mathcal E),
\]
and
\[
  0\le\frac{|\mu|(B(a,t))}{t}
  \le\int_{B(a,t)}\frac{d|\mu|(w)}{|w-a|}
  \xrightarrow{t\to0^+}0.
\]
If \(\Delta(a)\) is invertible at one such point, then
Lemma~\ref{model:regular} gives \(x\ne0\) and \(C_x<\infty\) such that for any $p\in\mathbb C[z]$,
\[
  \|p(T)x\|\le  C_x\|p\|_{\overline U}.
\]
Hence \(x\in X_T(\overline U)\setminus\{0\}\).

Suppose now that \(\Delta\) is singular \(\eta\)-almost everywhere.
By Lemma~\ref{model:regularized-transform}\textup{(1)},
\(\operatorname{ind}\Delta(z)=0\), so
\[
  \dim\ker\Delta(z)=\dim\ker\Delta(z)^*>0
  \qquad\text{for \(\eta\)-almost every }z.
\]
Let \(P(z)\) be the orthogonal projection onto \(\ker\Delta(z)\).
For each \(v\in\mathcal E\),
\[
  P(z)v=\lim\limits_{n\to\infty}
  \bigl(I+n\Delta(z)^*\Delta(z)\bigr)^{-1}v,
\]
and hence \(z\mapsto P(z)v\) is strongly measurable. Let \(\textcolor{black}{\{e_j\}}\) be an
orthonormal basis of \(\mathcal E\), then
\[
  \{z\in U:\ker\Delta(z)\ne\{0\}\}
  =\bigcup_{j\ge1}\{z\in U:P(z)e_j\ne0\}.
\]
Since \(\eta(U)>0\), choose \(v=e_j\) such that \(f(z)=P(z)v\)
satisfies
\[
  \eta\{f\ne0\}>0,\,
  \|f(z)\|\le1,\text{ and } \Delta(z)f(z)=0.
\]
For a bounded Bochner-measurable \(\mathcal E\)-valued field \(b\) on \(U\),
\textcolor{black}{set}
\[
  G_b^\eta(w)=\int_U\frac{b(z)}{z-w}\,d\eta(z),\text{ and }
  J_\eta b=\Gamma_A(G_b^\eta|_\Lambda).
\]
For \(r>0\), \textcolor{black}{set}
\(G_{b,r}^\eta(w)=\int_U k_r(z-w)b(z)\,d\eta(z)\).
Since \(k_r\) is bounded and continuous, \(G_{b,r}^\eta\) is
norm-continuous. By \(\eta\le\tau\) and Lemma~\ref{cau:energy},
\[
\begin{aligned}
  \|G_b^\eta-G_{b,r}^\eta\|_\infty
  \le \|b\|_\infty\sup_w
       \int_{U\cap B(w,r)}\frac{d\eta(z)}{|z-w|}
  \le b_\tau(r)\|b\|_\infty\xrightarrow{r\to0^+}0.
\end{aligned}
\]
Thus \(G_b^\eta\) is continuous in norm, and
\eqref{model:spectral-bound} yields
\[
  \|J_\eta b\|
  \le\|A\|_2\|G_b^\eta\|_\infty
  \le\|A\|_2 C_\tau\|b\|_\infty.
\]

Since
\[
  \int_\Lambda\int_U
  \frac{\|b(z)\|}{|z-w|}\,d\eta(z)\,d|\mu|(w)
  \le C_\tau\|b\|_\infty|\mu|(\mathbb C)<\infty,
\]
by Fubini's theorem and 
\eqref{model:spectral-identities} 
\[
  BJ_\eta b
  =\int_\Lambda d\mu(w)\,G_b^\eta(w)
  =\int_U\left(\int_\Lambda\frac{d\mu(w)}{z-w}\right)b(z)\,d\eta(z)
  =\int_U D(z)b(z)\,d\eta(z).
\]
Since \(G_{zb}^\eta(w)-wG_b^\eta(w)=\int_U b\,d\eta\), we obtain

 $$
  NJ_\eta b=J_\eta(zb)-A\int_U b\,d\eta,
$$
and 
$$
  TJ_\eta b=NJ_\eta b+ABJ_\eta b
            =J_\eta(zb)-A\int_U\Delta(z)b(z)\,d\eta(z).
$$
Set \(x=J_\eta f\). Since \(Df=f=P v\) almost everywhere,
\[
\begin{aligned}
  \langle Bx,v\rangle
  &=\int_U\langle D(z)f(z),v\rangle\,d\eta(z)\\
  &=\int_U\langle P(z)v,v\rangle\,d\eta(z)
   =\int_U\|P(z)v\|^2\,d\eta(z)>0.
\end{aligned}
\]
Thus \(x\ne0\). For every \(p\in\mathbb C[z]\), the identity
\(\Delta f=0\) gives
\[
  TJ_\eta(pf)
  =J_\eta(zpf)-A\int_U p(z)\Delta(z)f(z)\,d\eta(z)
  =J_\eta(zpf).
\]
Induction on monomials and linearity therefore yield
\(p(T)x=J_\eta(pf)\). Consequently,
\[
  \|p(T)x\|
  \le\|A\|_2C_\tau\|pf\|_\infty
  \le\|A\|_2C_\tau\|f\|_\infty\|p\|_{\overline U},
\]
so \(x\in X_T(\overline U)\setminus\{0\}\).

For the reflected measure \(\tau^\#(F)=\tau(F^\#)\), we have
\(C_{\tau^\#}=C_\tau\) and \(b_{\tau^\#}(r)=b_\tau(r)\).
Applying the same construction to \(T^*\) gives the
corresponding local vectors for \(T^*\). The disjoint-disk argument
in the proof of Proposition~\ref{eng:positive} now gives a closed
subspace \(\mathcal L\subset H\) satisfying
\[
  \dim\mathcal L=\dim\mathcal L^\perp=\infty,
\]
which is a hyperinvariant subspace of $T$.
\end{proof}

\begin{proof}[Proof of Theorem~\ref{main:half}]
\textcolor{black}{Let} \(\Lambda=\sigma(N)\). Since \(\dim_{\mathrm H}\Lambda>1\), choose
\(1<s<\dim_{\mathrm H}\Lambda\). By Frostman's lemma
\cite[Theorem~2.7]{Mattila}, there is a probability measure \(\omega\)
supported on \(\Lambda\) such that \(\omega(B(w,t))\le Ct^s\) for
\(w\in\mathbb C\) and \(0<t\le1\). Hence, by Tonelli's theorem,
\[
  I_1(\omega)
  =\int_\Lambda\int_0^\infty
    \frac{\omega(B(w,t))}{t^2}\,dt\,d\omega(w)
  \le C\int_0^1 t^{s-2}\,dt+\int_1^\infty t^{-2}\,dt<\infty.
\]
By Theorem~\ref{main:half1}, \(N+K\) has a hyperinvariant half-space.
\end{proof}

\medskip
\noindent\textbf{Declaration of AI Assistance.}
The authors used ChatGPT (OpenAI) to discuss and develop mathematical arguments and improve the writing. We carefully reviewed every suggestion given by AI and independently checked the mathematical content. Authors take full responsibility for the final manuscript.

\end{document}